\documentclass[12pt]{amsart}
\usepackage{amsmath,amssymb,graphics}

\DeclareMathOperator{\dist}{dist}
\begin{document}
\newtheorem{thm}{Theorem}[section]
\newtheorem{lem}[thm]{Lemma}
\newtheorem{dfn}[thm]{Definition}
\newtheorem{cor}[thm]{Corollary}
\newtheorem{conj}[thm]{Conjecture}
\newtheorem{clm}[thm]{Claim}
\theoremstyle{remark}
\newtheorem{exm}[thm]{Example}
\newtheorem{rem}[thm]{Remark}
\def\N{{\mathbb N}}
\def\Q{{\mathbb Q}}
\def\R{{\mathbb R}}
\def\C{{\mathbb C}}
\def\P{{\mathbb P}}
\def\Z{{\mathbb Z}}
\def\v{{\mathbf v}}
\def\r{{\mathbf r}}
\def\x{{\mathbf x}}
\def\O{{\mathcal O}}
\def\M{{\mathcal M}}
\def\tr{\mbox{Tr}}
\def\ord{{\mbox{ord}}}
\def\qed{{\tiny $\clubsuit$ \normalsize}}

\title[Rational approximations]
{A conjecture on rational approximations to rational points}
\author{David McKinnon}
\address{Department of Pure Mathematics \\
University of Waterloo \\
Waterloo, ON\ \  N2L 3G1 \\
Canada}
\email{dmckinnon@math.uwaterloo.ca}
\thanks{This research supported in part by a generous grant from NSERC}
\date{\today}

\begin{abstract}
In this paper, we examine how well a rational point $P$ on an
algebraic variety $X$ can be approximated by other rational points.
We conjecture that if $P$ lies on a rational curve, then the best
approximations to $P$ on $X$ can be chosen to lie along a rational
curve.  We prove this conjecture for a wide range of examples, and
for a great many more examples we deduce our conjecture from Vojta's
Main Conjecture.
\end{abstract}

\maketitle

\section{Introduction}

The distribution of rational points on an algebraic variety $X$ is
very subtle.  For example, it can often happen that there is a proper
subvariety $Y$ of $X$ such that the set of rational points of $Y$ has
density one in the set of rational points of $X$, where density is
used in the sense of Weil heights.  Roughly speaking, such
subvarieties $Y$ are called accumulating subvarieties.  Thus, if one
is interested in the arithmetic of $X$, one must first identify which
rational points of $X$ lie on $Y$, and which do not.

Unfortunately, even the purely geometric problem of identifying
potential accumulating subvarieties $Y$ can be very difficult.  It
would be helpful to have a local and arithmetic criterion to identify
points $P$ which lie on accumulating subvarieties.  In this paper, we
do not quite manage to construct such a criterion, but we do identify
an invariant, called the approximation constant of $P$ on $X$ with
respect to a divisor $D$, which describes how well $P$ can be
approximated by other rational points of $X$.  (See
Definition~\ref{bigdef} for details.)  If $P$ lies on an accumulating
subvariety, then it can be well approximated by other rational points,
and therefore will have a small approximation constant.

In the course of computing this constant in many examples, it became
clear that sequences which best approximate a rational point $P$ tend
to lie along curves.  While there are counterexamples to show that
this principle cannot hold in general, we are able to formulate a
conjecture which predicts that it should be true whenever $P$ lies on
a rational curve defined over the field of coefficients (see
Conjecture~\ref{ratcurve}).  We prove this conjecture for a wide
range of examples.

The basic technique we use for proving Conjecture~\ref{ratcurve} is
inductive.  We start by proving the conjecture for $\P^n$ (see
Theorem~\ref{ratapprox}).  We then prove a number of inductive
results which enable us to take advantage of the structure of the
N\'eron-Severi group of $X$.  Using these inductive results, we 
prove Conjecture~\ref{ratcurve} for a wide range of rational surfaces
by means of a careful analysis of the nef cone.

In the final section, we make several remarks about further directions
in which to attack Conjecture~\ref{ratcurve}, and drawbacks of our
current techniques.  We describe some cases which are unlikely to be
proven using our current techniques, and we further deduce
Conjecture~\ref{ratcurve} from Vojta's Main Conjecture (Conjecture
3.4.3 of \cite{Vo}) for a generic variety of non-negative Kodaira
dimension.

I am grateful to several people for helpful conversations about this
material, including Doug Park, Mike Roth, and Cam Stewart.  I would
particularly like to thank Kevin Hare for writing some extremely
useful computer programs, without which many of the theorems in this
paper would still be conjectures, and the anonymous referee, whose
invaluable suggestions improved the manuscript enormously.

\section{Lines}

In this paper, all heights are absolute and multiplicative.  We fix a
number field $k$.

\begin{thm}\label{ratapprox}
Let $P\in\P^n(k)$ be any $k$-rational point.  Let $\{P_i\}$ be any
sequence of $k$-rational points, and let $S$ be the set of archimedean
places of $k$, with the convention that pairs of complex conjugate
embeddings count as one place.  Suppose there is a positive real
constant $c\in\R$ such that for all $i$, we have:
\[\left(\sum_{v\in S}\dist_v(P,P_i)\right)H(P_i)\leq c\]
where the function $\dist_v$ denotes the local distance function on
$\P^n(k)$ induced by $v$.  Then the set $\{P_i\}$ is a subset of the
finite union of lines through $P$ of Pl\"ucker height at most $N(P)c$,
where $N(P)$ is a constant which does not depend on $c$.
\end{thm}

\noindent
{\it Proof:} \/ Without loss of generality, we may assume that $P=[0:0:\ldots
:0:1]$, and that none of the $P_i$ lie on the line $x_n=0$.  Fix a set $R$
of representatives of the class group of $k$.

Write $P_i = [a_{i1}:\ldots:a_{in}:b_i]$, where the $a_{ij}$ and $b_i$
are integers (that is, elements of $\O_k$) with
$\gcd(a_{i1},\ldots,a_{in},b_i)\in R$.  Let $v\in S$ be any
archimedean place of $k$.  Up to multiplication by a bounded function
(which won't affect the conclusion of the theorem), we may write the
$v$-distance as:
\[\dist_v(P,P_i) = \max_j\{|a_{ij}/b_i|_v\}\]
and the height as:
\[H(P_i) = \left(\prod_{v\in S}\max\{|a_{i1}|_v,\ldots,|a_{in}|_v,|b_i|_v\}
\right)^{1/[k:\Q]}\] 
Note that with our choice of representation, the non-archimedean
absolute values do not contribute more than a (multiplicatively)
bounded function to the height.

Now suppose that $\sum_{v\in S}\dist_v(P,P_i)H(P_i)\leq c$ for all
$i$.  Then we get:
\[\left(\sum_{v\in S}\max\{|a_{ij}/b_i|_v\}\right)^{[k:\Q]}
\prod_{v\in S}\max\{|a_{i1}|_v,\ldots,|a_{in}|_v,|b_i|_v\} \leq c^{[k:\Q]}\]
and hence
\[\left(\sum_{v\in S}\max\{|a_{ij}|_v\}\right)^{[k:\Q]}
\max\{|a_{i1}/b_i|_v,\ldots,|a_{in}/b_i|_v,1\} \leq c^{[k:\Q]}\]
and so {\it a fortiori}, we obtain:
\[\sum_{v\in S}\max\{|a_{ij}|_v\}\leq c\]
In particular, there are only finitely many choices for each $a_{ij}$,
and each choice corresponds to a line of height at most $c^{[k:\Q]}$.
\qed

\vspace{.1in}

In light of this theorem, we make the following definitions:

\begin{dfn}
Let $X\subset\P^n_k$ be an algebraic variety defined over $k$, and let
$S$ be the set of archimedean absolute values on $k$.  Define:
\[\dist(P,Q) = \sum_{v\in S} \dist_v(P,Q)\]
for any rational points $P$ and $Q$ on $X$.  Note that this is a well
defined distance function on the metric space $\prod_{v\in S} X(k_v)$,
where $k_v$ is the completion of $k$ with respect to $v$.
\end{dfn}

In this paper, all limits will be computed with respect to the topology
induced by the distance function $\dist$.

\begin{dfn}\label{appconst}
Let $X$ be an algebraic variety defined over $k$, and let $P\in X(k)$
be any rational point.  Let $D$ be any divisor on $X$, with
corresponding height function $H_D$.  Assume that there is a positive
constant $c$ such that $H_D(Q)>c$ for all $Q$ in some Zariski open
neighbourhood of $P$.  (This will be satisfied if, for example, $D$ is
ample, or more generally, if some multiple of $D$ is basepoint free.)
For any sequence $\{P_i\}\subset X(k)-\{P\}$ with $P_i\rightarrow P$,
define the {\em approximation constant on $X$ of $\{P_i\}$ with
respect to $P$ and $D$} to be the smallest non-negative real number
$\alpha$ such that:
\[\limsup_{i\rightarrow\infty}\dist(P,P_i)^\alpha H_D(P_i)<\infty\]
If there is no such smallest non-negative real number, then the sequence
does not have an approximation constant.
\end{dfn}

Notice that although the definition of $\dist$ depends on the choice of
embedding of $X$ in projective space, the definition of the approximation
constant does not.  This is because two different embeddings of $X$ are
diffeomorphic, and hence distances change by no more than a multiplicative
function bounded away from zero and infinity.  In this paper, this
ambiguity will never be significant, so we will refer to $\dist$ as a
function independent of the embedding of $X$.

\begin{exm}\label{noconst}
For example, let $P=[0\colon0\colon1]\in\P^2$, and let $f(n)=n/\log n
$.  On any $\Q$-rational line $L$ through $P$, we can find an infinite
sequence of rational points $\{P_i\}$ such that $\dist(P,P_i)
H(P_i)\leq c_L$, where $c_L$ is a constant that depends on $L$.  (See
for example Theorem~2.1 of \cite{M1}.)  Therefore, there is a positive
constant $C$ and an infinite sequence $\{Q_n\}\subset\P^2$ of rational
points such that $\dist(Q_n,P)f(H(Q_n))<C$ for all $n$.  Moreover, we
can choose the points $\{Q_n\}$ such that no two of them are collinear
with $P$.  Thus, by Theorem~\ref{ratapprox}, the approximation
constant of this sequence, if it exists, is greater than one.  But it
is clear that for every $\epsilon >0$, there exists a positive
constant $C_\epsilon$ such that
$\dist(Q_n,P)^{1+\epsilon}H(Q_n)<C_\epsilon$, so any constant of
approximation for this sequence is at most one.  We conclude that the
sequence $\{Q_n\}$ does not have a constant of approximation.
\end{exm}

\begin{dfn}\label{bigdef}
Let $X$, $P$, $D$, and $H_D$ be as in Definition~\ref{appconst}.
Define the {\em approximation constant of $P$ on $X$ with respect to
$D$} to be the minimum (if it is achieved) of all approximation
constants on $X$ of sequences $\{P_i\}\subset X(k)-\{P\}$ with respect
to $D$.  A {\em sequence of best approximation to $P$ with respect to
$D$} is a sequence whose corresponding approximation constant is equal
to the approximation constant of $P$.  A {\em curve of best
approximation to $P$} is a curve on $X$ passing through $P$ that contains
a sequence of best approximation to $P$.
\end{dfn}

Notice that by this definition, sequences with no approximation
constant (such as the one described in Example~\ref{noconst}) cannot
be a sequence of best approximation to a point $P$.  If
Conjecture~\ref{ratcurve} is true, then Theorem~\ref{curve} makes it
clear that this is a reasonable restriction.

We may therefore improve Theorem~\ref{ratapprox} as follows.

\begin{thm}\label{projlines}
Let $P\in\P^n(k)$ be any rational point.  Then there exists a sequence
of best approximation to $P$.  Moreover, any such sequence is a subset
of a finite set of lines and has constant of approximation equal to 1,
and any line through $P$ is a curve of best approximation to $P$.
\end{thm}

\noindent
{\it Proof:} \/ We first prove the theorem for $\P^1$.  Without loss
of generality, we may assume that $P=[0:1]$.  Fix a set of
representatives $R$ of the class group of $k$.  Let $Q=[a\colon b]$ be
any point in $\P^1(k)$, where $a$ and $b$ lie in $\O_k$ with $(a,b)\in
R$ and $a\neq 0$.  Since $\{[1/i:1]\}$ is clearly a sequence with
constant of approximation equal to 1, it suffices to show that
$\dist(Q,P)H(Q)\geq C_k$, where $C_k$ is a constant depending only on
$k$.

Since the arithmetic mean is always greater than the geometric mean,
for $Q$ close enough to $P$ we compute:
\[\dist(P,Q)=\sum_v|a/b|_v\geq [k:\Q]\prod_v|a/b|_v^{1/[k:\Q]}\]
\[=[k:\Q](N(a)/N(b))^{1/[k:\Q]}\]
where all sums and products are over archimedean places $v$, and $N$ denotes
the norm $N_{k/\Q}$.  We also compute:
\[H(Q)=(\prod_v\max\{|a|_v,|b|_v\})^{1/[k:\Q]}\geq N(b)^{1/[k:\Q]}\]
where the product again ranges over archimedean places $v$.  Note that
since $(a,b)\in R$, if we change the representation of $Q$, then we
will change the ideal $(a,b)$ by multiplication by an element of
$k^*$.  Since no two elements of $R$ represent the same ideal class,
it follows that the representation changes by multiplication by a unit
of $\O_k$, and so $H(Q)$ does not change.  Since $a\in\O_k$ is
nonzero, it follows that $N(a)\geq 1$, so that:
\[\dist(P,Q)H(Q)\geq [k:\Q]\]
and we have proven the result for $\P^1$.

In general, if $P\in\P^n$ is any rational point, we can find an
embedding $\phi\colon\P^1\hookrightarrow\P^n$ such that $\phi(\P^1)$
is a line, and $\phi(0:1)=P$.  Since $\phi$ only changes distances and
heights by multiplication by a function bounded away from zero and
infinity, it follows that a sequence of best approximation to $P$
along the line $\phi(\P^1)$ has constant 1.  But by
Theorem~\ref{ratapprox}, any sequence with approximation constant at
most 1 must lie in a finite set of lines through $P$.  The result
follows.  \qed

\vspace{.1in}

In light of Theorem~\ref{ratapprox}, we make the following conjecture:

\begin{conj}\label{ratcurve}
Let $V$ be an algebraic variety defined over $k$, and $D$ any ample
divisor on $V$.  Let $P$ be any $k$-rational point on $V$.  Assume
that there is some rational curve $C$ on $V$, defined over $k$, such
that $P\in C(k)$.  Then a sequence of best approximation to $P$ on $V$
with respect to $D$ exists, and may be chosen to lie along some
rational curve through $P$.
\end{conj}

It is easy to see that if we remove the hypothesis that $P$ lie on a
rational curve, then the conjecture is false.  Consider the case of a
simple abelian variety $V$ of dimension greater than one, and $P$ any
$k$-rational non-torsion point.  Then $P$ is certainly a limit point
of the set $V(k)$ in the archimedean topology, but any curve through
$P$ must have geometric genus at least two, which by Faltings' Theorem
can contain only finitely many $k$-rational points.

The following theorem is quite useful:

\begin{thm}\label{curve}
Let $C\subset\P^n$ be any irreducible rational curve of degree $d$
defined over $k$, and let $P\in C(k)$ be any $k$-rational point.  Let
$f\colon\P^1\rightarrow C$ be the normalization map, and let $m$ be
the maximum multiplicity of a branch of $f$ over $P$.  Then a sequence
of best approximation to $P$ on $C$ has constant of approximation
$d/m$.
\end{thm}

\noindent
{\it Proof:} \/ If $P$ is a smooth point of $C$, then the result
follows trivially from the observation that $f$ changes distances to
$P$ only by a multiplicative function bounded away from zero and
infinity, and $f$ raises heights to the power $d$ (up to a
multiplicative function bounded away from zero and infinity).

If $P$ is a singular point of $C$, then up to a multiplicative
function bounded away from zero and infinity, we have
$\dist(f(P),f(Q))=\dist(R,Q)$ near $P$, where $R$ is some point in the
finite set $f^{-1}(P)$.  Even though the point $R$ depends on the
point $Q$ near $P$, it nevertheless provides only a finite number of
alternatives for $\dist(f(P),f(Q))$, so that any sequence of best
approximation to $P$ must have an infinite subsequence which is a
sequence of best approximation to some element of $f^{-1}(P)$.

Thus, let $R\in f^{-1}(P)$ be any point, and let $m$ be the
multiplicity of the branch of $f$ through $R$.  Assume without loss of
generality that $R=0$.  Near $R$, in affine coordinates, the function
$f$ can be written as $f(x)=(f_1(x),\ldots,f_n(x))$, where
$C\subset\P^n$ and the $f_i$ are rational functions in $x$.  The
multiplicity $m_R$ of $f$ at $R$ is equal to
$\min_i\{\ord_0(f_i(x))\}$, where $\ord_0(f_i)$ denotes the largest
power of $x$ which divides the numerator of $f_i$ (in lowest terms).
Let $i$ be some index which achieves the minimum, and write
$f_i(x)=x^{m_R}+O(x^{m_R+1})$.  Then as $Q=x$ approaches $R=0$, the
distance between $Q$ and $R$ is $\sum_v|x|_v$, and the distance
between $f(Q)$ and $f(R)$ is (up to a bounded multiplicative constant)
$\sum_v|x^{m_R}|_v\ll\gg(\sum_v|x|_v)^{m_R}$.  The result now follows
from Theorem~\ref{ratapprox}.  \qed

\vspace{.1in}

Theorem~\ref{curve} explains why Conjecture~\ref{ratcurve} does not
refer specifically to rational curves of minimal $D$-degree.  A
rational point $P$ may be better approximated along a curve of higher
degree on which it is a cusp of high multiplicity than along a smooth
curve of much lower degree.

\begin{exm}\label{cusp}
Let $C$ be the plane curve $y^2z=x^3\subset\P^2$, and $P$ the cusp
$[0:0:1]$.  For $k=\Q$, it is not hard to see that a sequence of best
approximation to $P$ along $C$ is given by the sequence
$[1/B^2:1/B^3:1]$, which has constant of approximation $3/2$.  

However, let $f\colon\P^1\rightarrow C$ be the normalization map.
Then $f^*\O(1)$ is $D=\O(3)$.  A sequence of best $D$-approximation to
$P$ along $\P^1$ has constant $3$, not $3/2$, since $D$ corresponds to an
embedding of $\P^1$ in $\P^3$ as a twisted cubic curve.  This is
because in Definition~\ref{appconst}, we define $\dist$ in terms of an
embedding of $\P^1$, not an arbitrary morphism, whereas in this
example, we deal with the distance inherited from $\P^2$.

In other words, if $\tilde{W}\to W$ is the normalization of a
subvariety of a smooth variety $V$, then the function $\dist\colon
\tilde{W}\times \tilde{W}\rightarrow\R$ may not agree with the
pullback of the function $\dist\colon V\times V\rightarrow\R$ to
$\tilde{W}\times \tilde{W}$.  If $W$ is a curve, then the proof of
Theorem~\ref{curve} describes how these two functions differ, but if
$W$ has higher dimension, then one must be more careful.
\end{exm}

Despite the unpleasant possibility described by Example~\ref{cusp},
we can relate Conjecture~\ref{ratcurve} to the problem of finding
accumulating curves through a point $P$.

\begin{dfn}
Let $S\subseteq\P^n(k)$ be any set of $k$-rational points.  The
counting function for $S$ is defined to be:
\[N_S(B)=\#\{P\in S\mid H(P)\leq B\}\]
where $H$ denotes the standard height function on $\P^n$.
\end{dfn}

\begin{dfn}
Let $V\subseteq\P^n$ be an algebraic variety defined over $k$, and let
$W\subset V$ be any proper closed subset.  Then $W$ is said to be an
accumulating subvariety of $V$ if and only if:
\[N_{V-W}(B)=o(N_W(B))\]
where by $N_{V-W}$ and $N_W$ we mean $N_{V(k)-W(k)}$ and $N_{W(k)}$,
respectively.  Roughly speaking, an accumulating subvariety of $V$ is
a proper closed subset $W$ such that asymptotically, there are more
rational points of bounded height on $W$ than there are off $W$.
\end{dfn}

\begin{thm}
Let $V\subseteq\P^n$ be a smooth algebraic variety defined over $k$,
and let $P\in V(k)$ be any $k$-rational point.  Assume that
Conjecture~\ref{ratcurve} is true for $P$ on $V$ --- that is, assume
that there is a rational curve $C$ through $P$ which contains a
sequence of best approximation to $P$.  Assume further that $C$ has
only ordinary singularities at $P$.  If $P$ lies on an accumulating
curve of $V$, then that curve contains $C$ as an irreducible
component.
\end{thm}

\noindent
{\it Proof:} \/ First, note that since $V$ contains a rational curve
$C$ through $P$, it follows that any accumulating curve containing $P$
must be a rational curve, and moreover any component containing $P$
must be a rational curve of minimal degree through $P$.  By
Theorem~\ref{curve}, since $C$ has only ordinary singularities through
$P$, it is clear that any rational curve through $P$ must have degree
at least $\deg C$.  Thus, $C$ is a rational curve of minimal degree
through $P$, and therefore any accumulating curve through $P$ must
contain $C$ as an irreducible component.  \qed

\section{Other Varieties}

Next, we consider the question of rational approximation of rational
points on more general varieties.  We begin with a straightforward but
surprisingly useful result on products of varieties:

\begin{thm}\label{product}
Let $X$ and $Y$ be smooth algebraic varieties defined over a number
field $k$, and let $P$ and $Q$ be $k$-rational points on $X$ and $Y$,
respectively.  Let $L_X$ and $L_Y$ be divisors on $X$ and $Y$,
respectively, whose corresponding height functions $H_X$ and $H_Y$ are
bounded below by a positive constant in some Zariski open
neighbourhoods of $P$ and $Q$, respectively.  Let $\{P_i\}$
(respectively $\{Q_i\}$) be a sequence of best $L_X$-approximation
(respectively best $L_Y$-approximation) for $P$ (respectively $Q$),
with constant of approximation $\alpha$ (respectively $\beta$).  Then
either $\{(P_i,Q)\}$ or $\{(P,Q_i)\}$ is a sequence of best
$L$-approximation for $(P,Q)$ on $X\times Y$, where
$L=\pi_1^*L_X\otimes\pi_2^*L_Y$ ($\pi_i$ is the projection onto the
$i$th factor).  Moreover, any sequence $\{(S_i,T_i)\}$ whose constant
of approximation is less than $\alpha+\beta$ must have either $S_i=P$
or $T_i=Q$ for all but finitely many $i$.
\end{thm}

\noindent
{\it Proof:} \/ For any $k$-rational point $(S,T)$ on $X\times Y$, we
have $H_L(S,T) = H_X(S)H_Y(T)$, where $H_X$ and $H_Y$ are heights on
$X$ and $Y$ with respect to the divisors $L_X$ and $L_Y$,
respectively.  To prove the theorem, note that the non-negativity of
$\alpha$ and $\beta$ implies that it suffices to prove the last claim.

Let $\{(S_i,T_i)\}$ be any sequence of points in $(X\times
Y)(k)-\{(P,Q)\}$ with $(S_i,T_i)\rightarrow (P,Q)$.  We know that if
$\{S_i\}$ is infinite then $\dist(P,S_i)^{\alpha-\epsilon} H_X(S_i)$
is unbounded for any $\epsilon>0$.  Similarly, if $\{T_i\}$ is
infinite then $\dist(Q,T_i)^{\beta-\epsilon} H_Y(T_i)$ is unbounded
for any $\epsilon>0$.  Thus, if there are infinitely many $i$ for
which both $S_i\neq P$ and $T_i\neq Q$, then we deduce that the
following function is unbounded:
\[\max\{\dist(P,S_i),\dist(Q,T_i)\}^{\alpha+\beta-\epsilon}H_X(S_i)H_Y(T_i)\]
We may now deduce that the following function is unbounded:
\[\dist((P,Q),(S_i,T_i))^{\alpha+\beta-\epsilon}H_L(S_i,T_i)\]
which is to say that $\{(S_i,T_i)\}$ has constant of approximation 
at least $\alpha+\beta$.  \qed

\vspace{.1in}

We have the following immediate corollary:

\begin{cor}\label{divsum}
Let $X$ be a variety defined over $k$, and let $P\in X(k)$ be any
$k$-rational point.  Let $D_1$ and $D_2$ be two divisors on $V$ with
height functions $H_1$ and $H_2$ bounded below by a positive constant
in a neighbourhood of $P\in X(k)$.  Assume that $C_i$ is a curve of
best $D_i$-approximation to $P$ for each $i$, and let $D=a_1D_1+a_2D_2$
be any positive linear combination of $D_1$ and $D_2$.

\begin{itemize}
\item If $C_1=C_2=C$, then $C$ is also a curve of best
$D$-approximation to $P$.

\item If $C_1.D_2=0$, then either $C_1$ is a curve of best $D$-approximation
to $P$, or $C_2.D_1=0$ and $C_2$ is a curve of best $D$-approximation to 
$P$.
\end{itemize}

\end{cor}

The following theorem will be used in the proofs of Theorem~\ref{four}
and Theorem~\ref{onefibre}.

\begin{thm}\label{addeffective}
Let $X$ be a variety defined over $k$, and let $D$ and $E$ be divisors
on $X$.  Let $P\in X(k)$ be a rational point, and $C$ a curve of best
$D$-approximation to $P$.  Assume that $H_D$ and $H_{D+E}$ are bounded
below by positive constants in some neighbourhood of $P$, that $E$ is
effective, and that $C\cap E=\emptyset$.  Then $C$ contains a sequence
of best $(D+E)$-approximation to $P$.
\end{thm}

\noindent
{\it Proof:} \/ Since $C\cap E=\emptyset$, we can find a real,
positive constant $\alpha$ such that $H_{E}(Q)\in [1/\alpha,\alpha]$
for all $Q\in C(k)$.  Thus, the sequence of best $D$-approximation to
$P$ along $C$ will have the same constant of approximation with
respect to $D+E$ as with respect to $D$.  However, $\log H_E$ is
bounded from below away from $E$, and $P\not\in E(k)$.  It follows
that for any sequence $\{Q_i\}$ of rational points on $V$ with
$Q_i\rightarrow P$, if the quantity $\dist(P,Q_i)^\gamma H_D(Q_i)$ is
unbounded, then the quantity $\dist(P,Q_i)^\gamma H_{D+E}(Q_i)$ will
also be unbounded.  The theorem follows.  \qed

\vspace{.1in}

These theorems, surprisingly, give an immediate proof of
Conjecture~\ref{ratcurve} for any split, geometrically minimal surface
$X$ over $k$.  We say that a surface $X$ is split over $k$ if and only
if the inclusion of N\'eron-Severi groups $NS(X_k)\rightarrow
NS(X_\C)$ is an isomorphism -- that is, if every algebraic cycle on
$X_\C$ is numerically equivalent to some $k$-rational cycle.  The
proof relies crucially on the classification of minimal rational
surfaces over $\C$, which can be found in, for example, \cite{Be}.

\begin{cor}\label{minimal}
Conjecture~\ref{ratcurve} holds for any split, geometrically minimal
rational surface defined over $k$.  In particular, it holds for all
the Hirzebruch surfaces $H_n$ for $n\geq 0$.
\end{cor}

{\it Proof:} \/ The classification shows that any minimal rational
surface over $\C$ is either $\P^2$, $H_0=\P^1\times\P^1$, or a
Hirzebruch surface $H_n$ for some integer $n>0$.  The case of $\P^2$
is proven in Theorem~\ref{ratapprox}.  The case $\P^1\times\P^1$ is
immediate from Theorem~\ref{product}.  The case of $H_n$ follows from
Theorem~\ref{product} and Corollary~\ref{divsum} as follows.  Let
$A(H_n)$ be the closure in $N(H_n)=NS(H_n)\otimes\R$ of the cone of
ample divisors on $H_n$, where $NS(H_n)$ is the N\'eron-Severi group
of $H_n$.  Then $A(H_n)$ is generated by the two divisor classes
$D=S+nF$ and $F$, where $F=\pi^*(P)$ is a fibre of the map $\pi\colon
H_n\rightarrow\P^1$ and $S$ is the unique section of $\pi$ with
self-intersection $-n$.  The class $D$ is $f^*(L)$, where $f$ is the
contraction of $S$ to the vertex of the cone $C$ over a rational
normal curve of degree $n$ and $L$ is a line of the ruling of the
cone.

Consider the classes $D$, $F$, and $D+F$.  For any rational point $P$
on $H_n$, it is clear that a curve of best $F$-approximation to $P$ is
a fibre of $\pi$, and a curve of best $D$-approximation to $P$ is the
preimage of a line through $f(P)$ on the cone $C$ (if $P\not\in
S(k)$) or $S$ (if $P\in S(k)$).  If $P\not\in S(k)$, then these are
the same curve, so we are done by Corollary~\ref{divsum}.  If $P\in
S(k)$, then we simply note that the divisor $D+F$ corresponds to an
embedding of $H_n$ in projective space such that fibres of $\pi$ are
all lines, and $S$ is also a line, so they both contain sequences of
best $(D+F)$-approximation to $P$.  We again conclude by
Corollary~\ref{divsum}, since every very ample divisor on $H_n$ is a
positive linear combination either of $D$ and $D+F$ or $F$ and $D+F$.
\qed

\vspace{.1in}

\begin{thm}\label{simplefibres}
Let $n\geq 2$, and let $X$ be the blowup of the Hirzebruch surface
$H_n$ at $k<n$ points, no two of which lie in the same fibre of the
map to $\P^1$.  Conjecture~\ref{ratcurve} is true for $X$.
\end{thm}

\noindent
{\it Proof:} \/ First, it clearly suffices to assume that none of the
blown up points lies on the $(-n)$-section on $H_n$, since blowing up
points on this section will only increase $n$.  

Let $S$ be the class of the strict transform of the $(-n)$-section on
$H_n$, and let $F$ be the class of the total transform of a fibre of
the map from $H_n$ to $\P^1$.  For each of the $k$ reducible fibres, 
let $E_i$ and $F_i$ be the two components, of which $F_i$ is the one
which intersects $S$.

Let $\alpha\in\{0,1\}^k$ be any vector, and define a divisor
$D_\alpha$ on $X$ by:
\[D_\alpha = S+nF-\alpha\cdot(E_1,\ldots,E_k)\]
where the $\cdot$ denotes a formal dot product.

\begin{clm}
The effective cone of $X$ is generated by the divisors $E_i$, $F_i$,
and $S$.  The nef cone of $X$ is generated by the divisors $D_\alpha$
and $F$, where $\alpha$ ranges over all of $\{0,1\}^k$.
\end{clm}

\noindent
{\it Proof:} \/ The Picard rank of $X$ is $k+2$, so an arbitrary
divisor on $X$ can be written as $D=aS+bF+\sum f_iE_i$.  If $D$ is to
be ample, then its intersection with $S$, $E_i$, and $F_i$ must be 
positive, giving:
\[b-na>0\hspace*{1in}a+f_i>0\hspace*{1in}f_i<0\]
These inequalities define an open cone in $\R^{k+2}$ -- let $C$ be the
closure of this cone.  Then $C$ is finitely generated, and the
generating rays of $C$ are all intersections of $k+1$ of the 
hyperplanes $V=\{b-na=0\}$, $V_i=\{a+f_i=0\}$, and $W_i=\{f_i=0\}$.

Let $\v$ be a generator of $C$.  If $\v\not\in V$, then for some $j$,
$\v\in V_j\cap W_j$, since $\v$ must be contained in at least $k+1$
of the listed hyperplanes.  Thus, $\v$ is contained in the hyperplane
$a=0$.  Furthermore, since $\v$ is an extremal ray of $C$, there must
be some set of exactly $k+1$ of the $V_i$ and $W_i$ whose intersection
is the space generated by $\v$.  Since the vector $F$ is contained in 
the intersection of all the $V_i$ and $W_i$, it follows that $\v$ is
a positive multiple of $F$.

Thus, we may assume that $\v\in V$.  If $\v\in V_j\cap W_j$ for some $j$,
then $\v$ is again contained in $a=0$, and thus also in $b=0$.  But there
must also be some $\ell$ such that $\v\not\in V_\ell\cup W_\ell$, and
hence $f_\ell>0$ and $f_\ell<0$.  This is clearly a contradiction, and 
so for all $i$, $\v\not\in V_i\cap W_i$.

Let $\alpha\in\{0,1\}^k$ be the vector whose $i$th component is 0 if
$\v\in V_i$, and 1 if $\v\in W_i$.  Then $D_\alpha\in V_i$ if and only
if $\v\in V_i$, and similarly for $W_i$.  Since $D_\alpha\in V$, the
independence of the $k+1$ hyperplanes shows that $\v$ is a positive 
multiple of $D_\alpha$.

Thus, the two cones described in the claim are indeed dual to one
another.  To conclude the proof of the claim, it suffices to show that
every positive linear combination of the $D_\alpha$ and $F$ is ample.
First, it is clear that $F$ and $D_\alpha$ are basepoint free, since
$F$ corresponds to a morphism to $\P^1$, and $D_\alpha$ corresponds to
the morphism from $X$ to a cone, blowing down $S$ and exactly one
component of each reducible fibre.  We next note that $F^2=0$ and
$F.D_\alpha=1$, and:
\begin{eqnarray*}
D_\alpha.D_\beta & = & (S+nF)^2+(\alpha\cdot(E_1,\ldots,E_k))(\beta\cdot
(E_1,\ldots,E_k)) \\
& = & n - (\alpha\cdot\beta) \\
& \geq & n-k
\end{eqnarray*}
which is positive.  Therefore, by the Nakai-Moishezon criterion
(\cite{Ha}, Theorem V.1.10), any positive linear combination of $F$
and the $D_\alpha$ is ample, and thus the claim is proven.  \qed

\vspace{.1in}

We now prove Conjecture~\ref{ratcurve} for $X$.  If $P$ does not lie
on $S$ or any reducible fibre, then for each $\alpha$, a sequence of
best $D_\alpha$-approximation to $P$ is clearly contained in the
component $C$ of $F$, and certainly a sequence of best
$F$-approximation is also contained in $C$.  Thus, for any ample $D$,
$C$ is a curve of best $D$-approximation to $P$.  

If $P$ lies on exactly one $E_i$ or $F_i$ (and not $S$), then the same
analysis shows that this same $E_i$ or $F_i$ is a curve of best
$D$-approximation to $P$ for any ample $D$.  If $P$ lies on $S$ but no
$E_i$ or $F_i$, then notice that $S.D_\alpha=0$ for any $\alpha$, so
we conclude that $S$ is a curve of best $D_\alpha$-approximation to
$P$ for all $\alpha$.  Since $S.F=1$, and since $D_\alpha+F$ contracts
no curves through $P$, we conclude from Theorem~\ref{projlines} that
$S$ and $C$ are both curves of best $(D_\alpha+F)$-approximation to
$P$, for any $\alpha$.  Since any ample divisor $S$ is either a
positive linear combination of elements of the set $\{D_\alpha,
F+D_\alpha\}$, or a positive linear combination of elements of the set
$\{F,F+D_\alpha\}$, we conclude from Corollary~\ref{divsum} that
either $S$ or $C$ is a curve of best $D$-approximation to $P$ for any
ample $D$.

If $P=S\cap E_i$ for some $i$, then a similar analysis shows that for
any ample divisor $D$, either $S$ or $E_i$ is a curve of best
$D$-approximation to $P$.  In particular, if we divide the nef cone of
$X$ into two subcones according to which of $S$ or $E_i$ has smaller
degree, then it is straightforward to check that the curve of smaller
degree has degree zero or one with respect to each generator of the
corresponding subcone, and so by Theorem~\ref{projlines} and
Corollary~\ref{divsum}, we conclude that the curve of smaller
$D$-degree is always a curve of best $D$-approximation to $P$.  If
$P=E_i\cap F_i$ for some $i$, a similar calculation shows that either
$E_i$ or $F_i$ --- whichever has smaller $D$-degree --- is always a
curve of best $D$-approximation to $P$.  \qed

\vspace{.1in}

Theorem~\ref{simplefibres} can be used to show that any $k$-split
rational surface of Picard rank at most three satisfies
Conjecture~\ref{ratcurve}.

\begin{thm}\label{three}
Let $X$ be a $k$-split rational surface, of Picard rank at most three.  Then
Conjecture~\ref{ratcurve} is true for $X$.
\end{thm}

\noindent
{\it Proof:} \/ If the Picard rank is one, then $X$ is isomorphic to
$\P^2$.  If the Picard rank is two, then either it's geometrically
minimal, or else it's $\P^2$ blown up at a point $P$.  In this case,
however, $X=H_1$, and the result follows from Corollary~\ref{minimal}.

If the Picard rank of $X$ is three, then it must be the blowup of some
$H_n$ at some point $P$, with $n\geq 0$.  If $n=0$, then since any blowup
of $H_0=\P^1\times\P^1$ is also a blowup of $H_1$, we may instead take
$n=1$.  If $n>1$, then we can apply Theorem~\ref{simplefibres} directly.
Thus, we assume that $n=1$, and therefore that $X$ is the blowup of
$\P^2$ at two different points.

We begin by finding generators for the effective cone of $X$.  Let
$\pi\colon X\rightarrow\P^2$ be the blowing down map, and let $E_1$
and $E_2$ be the two exceptional divisors.  Let $S$ be the strict
transform of the line in $\P^2$ which joins the two blown up points.
If we write $L=S+E_1+E_2$, then $L$ is the class of $\pi^*\O(1)$.  Then
$L-E_1$ and $L-E_2$ correspond to morphisms $\psi_1$ and $\psi_2$ to
$\P^1$.  We have the following well known description of the effective
and nef cones of $X$.

\begin{clm}
The curves $S$, $E_1$, and $E_2$ generate the closed cone
$\overline{NE}(X)$ of effective divisors on $X$, and the classes $L$,
$L-E_1$, and $L-E_2$ generate the nef cone of $X$.
\end{clm}

We now prove Conjecture~\ref{ratcurve} for $X$.  Let $P\in X(k)$ be
any point, and let $D$ be any ample divisor.  If $P\neq S\cap E_i$,
let $C_i$ be the unique irreducible component of the fibre of
$\psi_i\colon X\to\P^1$ through $P$.  Each of $L$, $L-E_i$, and
$A=(L-E_1)+(L-E_2)$ either contracts $C_i$, or else contracts no
curves through $P$ and maps $C_i$ to a line.  (Note that if $P$ lies
on $S$, then $C_i=S$.)  Thus, by Theorem~\ref{projlines} and
Corollary~\ref{divsum}, we conclude that $C_i$ is a curve of best
$D$-approximation to $P$ for any $D$ in the cone generated by $L$,
$L-E_i$, and $A$.

This covers all cases except $P=S\cap E_i$ and $D$ is a positive
linear combination of $L$, $L-E_i$, and $A$.  In this case, notice
that by Theorem~\ref{product}, both $S$ and $E_i$ are curves of best
approximation to $P$ with respect to $2L-E_i$ and $L+A$.  Since $L$
contracts $E_i$, this means that if $D$ lies in the cone generated by
$L$, $2L-E_i$, and $L+A$, then $E_i$ is a curve of best
$D$-approximation to $P$ by Corollary~\ref{divsum}.  And since $L-E_i$
and $A$ both contract $S$, we conclude again by Corollary~\ref{divsum}
that $S$ is a curve of best $D$-approximation to $P$ if $D$ is a
positive linear combination of $L-E_i$, $2L-E_i$, $A$, and $L+A$.  This
concludes the proof of Conjecture~\ref{ratcurve} for $X$.  \qed

\vspace{.1in}

Proceeding to the other extreme from Theorem~\ref{simplefibres}, we
now prove a theorem about the case in which the map to $\P^1$ has
exactly one reducible fibre.

\begin{thm}\label{onefibre}
Let $X$ be a smooth rational surface obtained by a succession of
blowups of the Hirzebruch surface $H_n$, and let $\pi\colon
X\rightarrow\P^1$ be the associated map.  Assume that $\pi$ has only
one reducible fibre, with $m$ components, and assume $m<n$.  If every
multiple component of the reducible fibre intersects at least two
other components then Conjecture~\ref{ratcurve} is true for every
point of $X$.
\end{thm}

\noindent
{\it Proof:} \/ If any of the blown-up points of $H_n$ lie on the
$(-n)$-section $S$, then the strict transform of $S$ will have
strictly smaller self-intersection than $-n$, and $X$ can be obtained
as a blowup of $H_r$ for some $r>n$.  Thus, if we choose $n$ large
enough, we may therefore assume that none of the blown-up points of
$H_n$ lie on $S$.

Let us establish some notation:
\begin{itemize}
\item The components of the reducible fibre of $\pi$ are denoted by $E_1,
\ldots,E_m$, and we write $E_0=S$ for the unique $(-n)$-section of $\pi$.

\item The classes $E_0,\ldots,E_m$ are a basis for the N\'eron-Severi
group $NS(X)$, and we denote by $D_0,\ldots,D_m$ the dual basis with 
respect to the intersection pairing.  That is, $D_i.E_j=\delta_{ij}$.

\item We write $F=m_1E_1+\ldots +m_mE_m$, so that $m_i$ is the
multiplicity of $E_i$ as a component of the reducible fibre of $\pi$.
\end{itemize}

We proceed with a technical lemma.

\begin{lem}\label{inductivestep}
Assume without loss of generality that $E_m^2=-1$ and $E_m.S=0$, and
let $f\colon X\rightarrow Y$ be the map which blows down $E_m$.  Let
$\phi\colon Y\rightarrow\P^1$ be the map satisfying $\pi=\phi\circ f$,
let $E'_0$ be the $(-n)$-section of $\phi$, and let
$E_1',\ldots,E_{m-1}'$ be the components of the unique reducible
section of $\phi$, ordered so that $E_i'=f_*(E_i)$.

Let $\{D_0',\ldots,D_{m-1}'\}$ be the dual basis to
$\{E'_0,E_1',\ldots,E_{m-1}'\}$.  Then $D_i=f^*(D_i')$ if $0\leq i<m$.

Furthermore, $E_m$ can intersect either one or two other components
$E_i$.  If $E_m$ intersects one component $E_i$, then
$D_m=f^*(D_i')-E_m$.  If $E_m$ intersects $E_i$ and $E_j$, then
$D_m=f^*(D_i')+f^*(D_j')-E_m$.
\end{lem}

\noindent
{\it Proof of Lemma:} \/ The fact that $E_m$ can intersect at most two
components $E_i$ follows from the fact that no point of $Y$ can lie on
more than two components of the reducible fibre.  The $D_i$ are
defined by the property that $D_i.E_j=\delta_{ij}$, where we define
$E_0=S$.  It is a straightforward matter to verify this equality for
all $i$ and $j$.  \qed

\vspace{.1in}

The next step is to compute the effective cone of $X$.

\begin{lem}\label{effectcone}
Assume that for all $i$, some multiple of the divisor $D_i$ is
basepoint free.  Then the effective cone of $X$ is generated by the
components $\{E_1,\ldots,E_m\}$ of the reducible fibre and the
$(-n)$-section $S$ of $\pi$.  Furthermore, we have $D_0=F$ and
$D_1=S+nF$.
\end{lem}

\noindent
{\it Proof of Lemma:} \/ Let $E_1,\ldots,E_m$ be the set of components
of the reducible fibre, where $E_1$ is the unique component which
intersects $S$.  Let $A=(a_{ij})_{i,j=0}^m$ be the intersection matrix
of the $S$ and $E_i$; that is, let $a_{ij}=E_i.E_j$, where $E_0=S$.
The only entry of $A$ which is not $0$ or $\pm 1$ is $a_{00}=-n$.  We
therefore may regard $A$ as a specialization of the matrix $A_N$ with
entries in $\Z[N]$, where $(A_N)_{ij}=a_{ij}$ except that
$(A_N)_{00}=-N$.  Specializing $N=n$ reduces $A_N$ to $A$.

To check that $\{S,E_i\}$ generates the effective cone, it suffices to
check that the dual cone is the nef cone.  Since $\{S,E_i\}$ form a
basis for the vector space $NS(X)\otimes\R$, the dual cone will be
generated by the vectors $D_j = \sum b_{ij}E_i$, where
$B_N=(b_{ij}(N))_{i,j=0}^m$ is the inverse of the matrix $A_N$,
$b_{ij}=b_{ij}(n)$, and $E_0=S$.  To show that this dual is the ample
cone, it suffices to show that every positive linear combination of
the $D_j$ is ample.  By the Nakai-Moishezon criterion, since we assume
that some multiple of each $D_i$ is basepoint free, it suffices to
show that $b_{ij}=D_i.D_j\geq 0$ for all $i$ and $j$, and that for
each $i$, there is some $j$ for which $D_i.D_j>0$.

First, note that $D_0=F$, the class of a fibre of $\pi$.  To see this,
note that $D_0.E_i=0$ for all $i>0$, and therefore $D_0=\lambda F$ for
some $\lambda$ by the non-degeneracy of the intersection pairing.  Since
$D_0.S=1$, we conclude that $D_0=F$.  In particular, the first row and 
column of $B$ are independent of $N$, and $b_{0j}>0$ for all $j>0$.

Now let $i>0$.  Since the only entry of $A_N$ depending on $N$ is
$(A_N)_{00}=-N$, Cramer's Rule implies that the $i$th row $\r_i$ of
$A_N$ is $\r_i=N\r_i^\prime+\r_i^\circ$, where $\r_i^\prime$ and
$\r_i^\circ$ are independent of $N$.  Specializing to $N=n$, we obtain
$D_i=nD_i^\prime+D_i^\circ$.  (Note that $F.E_j=0$ for all $j>0$, so
that the $(0,0)$-cofactor matrix of $A_N$ has nontrivial kernel, and
therefore zero determinant.  We therefore conclude that the
determinant of $A_N$ is $\pm1$, and in particular independent of $N$.)
Furthermore, we know that if $\r^*_j$ denotes that $j$th column of
$A_N^{-1}$, then $\r_i\cdot\r^*_j=\pm D_i.E_j=\delta_{ij}$ is
independent of $N$ for all $j$, and hence $D_i^\prime.E_j=0$ for all
$j>0$.  Thus, $D_i^\prime=\lambda_i F$ for some $\lambda_i$.

Clearly $D_1=S+nF$, since $(S+nF).S=0$, $(S+nF).E_j=0$ if $j>1$, and
$(S+nF).E_1=1$.  Thus, $\lambda_1=1>0$.  Furthermore, since
$S+nF=D_1=\sum_j b_{j1}E_j$, this means that $b_{1j}=b_{j1}=\lambda_j
n$ is just $n$ times the multiplicity of $E_j$ in the fibre $F$ --- in
particular, $\lambda_j$ is a positive integer.  It therefore remains
only to show that the coefficient of $E_j$ in $D_i^\circ$ is at least
$-m\lambda_j$, and that at least one coefficient is strictly greater
than $-m\lambda_j$.  (Recall that $m$ is the number of components in
the reducible fibre.)

We proceed by induction on $m$.  If $m\leq 1$, the result follows
immediately from the preceding calculations, since $D_0$ and $D_1$
will be the full list of $D_i$'s.  For general $m$, note that there
must be some $E_i$ with $E_i^2=-1$ and $i\neq 1$ --- without loss of
generality, we may assume that $i=m$.  Let $f\colon X\rightarrow Y$ be
the map that blows down $E_m$.  Then $Y$ is a smooth rational surface
which admits a map $\phi \colon Y\rightarrow\P^1$ such that $\phi\circ
f=\pi$, and the unique singular fibre of $\phi$ has $m-1$ components.

Let $E^\prime_i=f_*(E_i)$ for $i<m$, and let $S'=f_*(S)$.  Denote the
dual basis of $\{S',E_i^\prime\}$ by $\{D_i^\prime\}$.  By induction,
the coefficient of $E^\prime_j$ in $D_i^{\prime\circ}$ is at least
$-(m-1)\lambda_j$, since the multiplicity of $E^\prime_j$ in $f_*(F)$
is also $\lambda_j$.  If $i<m$, then by Lemma~\ref{inductivestep}, the
coefficient of $E_j$ in $D_i^\circ$ is also at least
$-(m-1)\lambda_j\geq -m\lambda_j$.  If $i=m$, then the coefficient of
$E_j$ in $D_i$ is at least $-(m-1)\lambda_j-1\geq -m\lambda_j$.

Thus, if $m\leq n$, then $D_i.D_j\geq 0$ for all $i$ and $j$.  Since
the pairing is non-degenerate, for each $i$ there is some $j$ for
which $D_i.D_j>0$.  This concludes the proof of the lemma.  \qed

\vspace{.1in}

\begin{rem}\label{strictremark}
Note that the proof of Lemma~\ref{effectcone} also shows that
$D_i.D_j>0$ for all $i$ and $j$ except $i=j=0$.
\end{rem}

The next step is to prove that for each $i$, some multiple of $D_i$ is
basepoint free.  Together with Lemma~\ref{effectcone}, this will
complete the computation of the effective cone of $X$.  

Define a graph $G=(V,E)$ with vertex set $V=\{E_1,\ldots,E_m\}$ and
edge set $E=\{(E_i,E_j)\mid E_i.E_j>0\}$.  Since $E_i^2<0$, this
defines a simple graph, and moreover, it is a tree.  Root this tree at
$E_1$; that is, define a partial order on $V$ by $E_j\succeq E_i$ if
and only if there is a simple path (that is, a path with no repeated
vertices) from $E_1$ to $E_j$ which contains the vertex $E_i$.  Call
$E_i$ a leaf of $G$ if and only if $E_j\succeq E_i$ implies $i=j$.
The hypothesis of the theorem is equivalent to demanding that the
leaves of $G$ have multiplicity one in the fibre $F$.

\begin{lem}\label{leafbasepointfree}
Let $E_i$ be a leaf of $G$.  Then $D_i$ is basepoint free.
\end{lem}

\noindent
{\it Proof:} \/ We proceed by double induction on $n$ and $m$.  If
$n=1$ there there is nothing to prove.  If $m=1$, then $X$ is the
Hirzebruch surface $H_n$ and the lemma is clear.  If $n=2$, then we
must have $m\leq 1$, and the result is again clear.  Note also that if
$m>1$, then the sum of the multiplicities of the $E_i$ with $E_i^2=-1$
is always at least two, and that this property will be preserved by
blowing up.

Now consider a general $n>2$, and let $E_i$ be a leaf of $G$.  Since
$E_i$ has multiplicity one by hypothesis, there must be some other
$(-1)$-curve $E_j$.  (It is possible that $j=1$.)  Let $f\colon
X\rightarrow Y$ be the map that blows down $E_j$.  Then $Y$ will also
be a rational surface with a map $\phi\colon Y\rightarrow\P^1$ with at
most one reducible fibre.  The unique $(-\ell)$-section of $\phi$ will
satisfy $\ell\leq n$, and the reducible fibre will have $m-1$
components.  Therefore, by induction, the dual divisor $D_i^\prime$ to
$f_*(E_i)$ is basepoint free.  But $D_i=f^*(D_i^\prime)$, so $D_i$ is
basepoint free as well.  \qed

\vspace{.1in}

\begin{lem}\label{multiplegens}
Assume that $E_i\succeq E_j$.  Then $m_iD_j-m_jD_i$ is an effective
divisor supported on the reducible fibre of $\pi$.  Furthermore, if
$E_i$ and $E_j$ intersect, then $m_iD_j-m_jD_i=\sum_{E_t\succeq
E_i}m_tE_t$.
\end{lem}

\noindent
{\it Proof:} \/ Since $F.(m_iD_j-m_jD_i)=0$, it follows that
$m_iD_j-m_jD_i$ is supported on the reducible fibre of $\pi$.  Thus,
it remains only to show that it is effective.  It suffices to prove
the result in the case that $E_i.E_j=1$; a simple induction will prove
the general case from there.  Furthermore, if we blow down a
$(-1)$-curve $E_\ell$ on $X$ to obtain $Y$, then by
Lemma~\ref{inductivestep}, the truth of the lemma for $E_i$ and $E_j$
on $X$ will be equivalent to the truth of the lemma for the images of
$E_i$ and $E_j$ on $Y$, unless $i$ or $j$ equals $\ell$.  Thus, by
induction on the number $m$ of components of the reducible fibre, we
may assume that $E_i$ or $E_j$ is a $(-1)$-curve on $X$, that any
$(-1)$-curve in the reducible fibre is either $E_i$ or $E_j$, and that
$i\neq 1$ (because $E_i\preceq E_j$ and $i\neq j$).

Say $E_i$ is a $(-1)$-curve of multiplicity one.  Then it must be a
leaf of $G$, so $E_j$ must be the unique other component which
intersects $E_i$.  In that case, $E_j$ must also have multiplicity
one, so let $Y$ be the surface obtained by blowing down $E_i$.  By
applying Lemma~\ref{inductivestep}, we see that $D_j-D_i =E_i$, and
the lemma follows.

Now assume that $E_i$ is a $(-1)$-curve of multiplicity $m_i>1$.  Then
by hypothesis $E_i$ cannot be a leaf, so it must intersect two curves
$E_j$ and $E_\ell$, with $m_i=m_j+m_\ell$.  We have $E_\ell\succeq
E_i\succeq E_j$, so let $f\colon X\rightarrow Y$ be the map that blows
down $E_i$.  Then $Y$ is also a rational surface with $(-n)$-section
and a single reducible fibre with $m-1$ components, so
$\{S',E_1',\ldots,E_{m-1}'\}$ are generators of the effective cone of
$Y$, with $S'=f_*(S)$ and $E_r'=f_*(E_r)$ for any $r$.  Let
$\{D_0',\ldots,D_{m-1}'\}$ be the dual basis of
$\{S',E'_1,\ldots,E'_{m-1}\}$ in $\mbox{NS}(Y)$.  By
Lemma~\ref{inductivestep}, we have $D_i=D_j+D_\ell-E_i$.  By induction
on $m$, we have:
\[m_\ell D_j' - m_j D_\ell'=\sum_{E_t\succeq E_\ell} m_tE_t'\]
since the multiplicity of $E_r$ is the same as the multiplicity of
$E_r'$, for any $r\neq i$.  By Lemma~\ref{inductivestep} again,
pulling back to $X$ via $f$ gives:
\[m_\ell D_j - m_j D_\ell =(\sum_{E_t\succeq E_\ell} m_tE_t)+m_\ell E_i\]
and therefore:
\begin{eqnarray*}
m_iD_j-m_jD_i & = & m_iD_j-m_jD_j-m_jD_\ell +m_jE_i \\
& = & m_\ell D_j-m_jD_\ell + m_jE_i \\
& = & (\sum_{E_t\succeq E_\ell} m_tE_t) + m_iE_i \\
& = & \sum_{E_t\succeq E_i} m_tE_t
\end{eqnarray*}
as desired.  

Finally, we treat the case that $E_j$ is the only $(-1)$-curve in the
reducible fibre.  It is not a leaf, since $E_i\succeq E_j$.  If $j=1$,
then since $E_1=E_j$ has multiplicity one, then since $m>1$ ($i\neq
j$), there must be at least one other $(-1)$-curve $E_t$ on $X$.
Therefore $j\neq 1$, and let $E_i$ and $E_\ell$ be the two components
which intersect $E_j$.  (Note that $j\neq 1$ because $E_i\preceq
E_j$.)  We have $E_i\succeq E_j\succeq E_\ell$, $D_j=D_i+D_\ell -E_j$,
and $m_j=m_i+m_\ell$.  By induction and Lemma~\ref{inductivestep}, we
compute:
\[m_iD_\ell - m_\ell D_i = (\sum_{E_t\succeq E_i} m_tE_t)+m_iE_j\]
and therefore:
\begin{eqnarray*}
m_iD_j-m_jD_i & = & m_iD_i+m_iD_\ell-m_iE_j-m_jD_i \\
& = & (m_iD_\ell - m_\ell D_i) - m_iE_j \\
& = & \sum_{E_t\succeq E_i} m_tE_t
\end{eqnarray*}
and the lemma is proven.  \qed

\vspace{.1in}

\begin{lem}\label{basepointfree}
For any $i$, some positive multiple of $D_i$ is basepoint free.
\end{lem}

\noindent
{\it Proof of lemma:} \/ We have already proven this in the case that
$i=0$ (proof of Lemma~\ref{effectcone}) or $E_i$ is a leaf (see
Lemma~\ref{leafbasepointfree}).  By Lemma~\ref{multiplegens}, we know
that some multiple of each $D_i$ is effective, and in particular, the
base locus of $D_i$ is supported on the set $\cup_{E_t\succ E_i} E_t$.
It therefore suffices to find an effective divisor $E$, linearly
equivalent to some multiple of $D_i$, supported on the set
$S\cup(\bigcup_{E_t\not\succeq E_i} E_t)$.  We claim that for all $i$,
there is an effective divisor $C_i$ and a positive rational number
$\alpha_i$ such that $D_i=(m_iD_1-\alpha_iF)+C_i$, $C_i$ is supported
on the set $\cup_{E_t\not\succeq E_i} E_t$, and the base locus of a
suitable multiple of $m_iD_1-\alpha_iF$ is a subset of $S$.  In
particular, we will show that $\alpha_i/m_i\leq\ell_i$, where $\ell_i$
denotes the number of edges in the shortest path in the graph $G$ from
$E_1$ to $E_i$.

If $i=1$, this is trivial.  We proceed by induction on $\ell_i$, the
number of edges in the shortest path from $E_1$ to $E_i$ in the rooted
tree.  Let $E_j$ be the parent of $E_i$ --- that is, let $E_j$ be the
unique component with $E_j\prec E_i$ and $E_j.E_i=1$.  Then by
induction we can write $D_j=(m_jD_1-\alpha_jF)+C_j$, where $C_j$ is
effective and supported on the set $\cup_{E_t\not\succeq E_j} E_t$,
and the base locus of a suitable multiple of $m_jD_1-\alpha_jF$ is a
subset of $S$.  By Lemma~\ref{multiplegens}, we can write:
\[m_iD_j-m_jD_i=\sum_{E_t\succeq E_i} m_tE_t\]
We therefore deduce:
\begin{eqnarray*}
D_i & = & (1/m_j)(m_iD_j-\sum_{E_t\succeq E_i}m_tE_t) \\
& = & (1/m_j)(m_im_jD_1-m_i\alpha_jF+m_iC_j-\sum_{E_t\succeq E_i}m_tE_t) \\
& = & m_iD_1-(m_i\alpha_j/m_j+1)F+m_iC_j+\sum_{E_t\not\succeq E_i}m_tE_t
\end{eqnarray*}
Set $\alpha_i=m_i\alpha_j/m_j+1$.  Since $\alpha_j/m_j\leq\ell_j$, it
follows that $\alpha_i/m_i=\alpha_j/m_j+1/m_i \leq\ell_i=\ell_j+1$.
We have $\ell_i\leq m<n$ for all $i$, so a suitable multiple of
$D_1-\alpha_i/m_iF=S+(n-\alpha_i/m_i)F$ will be basepoint free away from $S$,
as desired.  This means that for each $i$, the base locus of a
suitable multiple of $D_i$ is contained in the union of $S$ and
$\cup_{E_t\not\succeq E_i}E_t$.  But we have already proven that the
base locus of a suitable multiple of $D_i$ is contained in the union
$\cup_{E_t\succ E_i} E_t$.  Since these two sets are disjoint, we
conclude that a large enough multiple of $D_i$ is basepoint free, as
desired.  \qed

\vspace{.1in}

We now conclude the proof of Theorem~\ref{onefibre}.  Choose any $P\in
X(k)$.  

\begin{clm}\label{vertexproof}
For any $i$, a sequence of best $D_i$-approximation to $P$ can be
found on the fibre of $\pi$ through $P$ (this includes possibly along
some $E_j$), or else on $S$.
\end{clm}

\noindent
{\it Solution:} \/ If $i=0$, then $D_i=F$, so the claim is trivial.
If $i=1$, then $D_i=S+nF$, and against the claim is clearly true.  For
general $i$, if $P$ lies on a curve $C$ with $C.D_i=0$, then clearly
$C$ will be a curve of best $D_i$-approximation to $P$.  However, such
a curve $C$ must be linearly equivalent to a nonnegative linear
combination of $\{E_j\mid j\neq i\}$, where as usual $E_0=S$.  If $C$
is not $S$ or some $E_i$, then $C$ must be nef, and so must be a
positive linear combination of the $D_i$.  Since $m<n$,
Remark~\ref{strictremark} shows that this is impossible.  Thus, if
$C.D_i=0$, then $C$ is either $S$ or some $E_i$ for $i\geq 1$.

Assume that $E_i$ is a leaf with $i>1$.  Then $D_i$ is basepoint free,
and $D_i.F=1$, so if $\phi_i\colon X\to\P^n$ is the morphism
associated to $D_i$, then $\phi_i$ maps the fibre of $\pi$ through $P$
to a line.  By Theorem~\ref{projlines}, this means that a sequence of
best $D_i$-approximation to $P$ either lies along the fibre of $\pi$
through $P$, or else along some irreducible curve $C$ with $C.D_i=0$,
which as noted must either be $S$ or some $E_i$.

If $E_i$ is not a leaf, then there is some $j$ for which $E_j$ is a
leaf and $E_j\succeq E_i$.  By Lemma~\ref{multiplegens}, we can write
$D_i=m_i D_j+\sum_{E_t\succeq E_j} m_tE_t$.  If $P$ does not lie on
any $E_t$ with $E_t\succeq E_j$, then Theorem~\ref{addeffective}
implies that either the fibre of $\pi$ through $P$ or $S$ is a curve
of best $D_i$-approximation to $P$.  If $P$ lies on some $E_t$ for
$t\neq i$, then $E_t$ is a curve of best $D_i$-approximation to $P$
because $E_t.D_i=0$ and some multiple of $D_i$ is basepoint free.  If
$P$ lies on $E_i$ but not on any other $E_t$, then as in the proof of
Lemma~\ref{basepointfree}, we can write $D_i=\alpha F+E$, where
$\alpha$ is a positive rational number and $E$ is an effective divisor
whose support does not include $E_i$.  Since $E_i$ is a curve of best
$F$-approximation to $P$, Lemma~\ref{multiplegens} implies that $E_i$
is also a curve of best $D_i$-approximation to $P$.  In all cases,
the claim is proven.  \qed

\vspace{.1in}

The claim, together with Corollary~\ref{divsum}, shows that if $D$ is
any ample divisor on $X$, then $D$ is a nonnegative linear combination
of the $D_i$, and so a sequence of best $D$-approximation to $P$ can
be found on the fibre of $\pi$ which contains $P$, unless $P$ lies on
$S$.  If $P$ lies on $S$, then $S$ is a curve of best $D_i$-approximation
to $P$ for all $i$, and so $S$ is also a curve of best $D$-approximation 
to $P$.  This concludes the proof of Theorem~\ref{onefibre}.  \qed

\vspace{.1in}

These same techniques will also prove Conjecture~\ref{ratcurve} for 
split rational surfaces of Picard rank four.

\begin{thm}\label{four}
Conjecture~\ref{ratcurve} is true for any split rational surface $X$
of Picard rank four.
\end{thm}

\noindent
{\it Proof:} \/ Every such rational surface is a blowup of a Hirzebruch
surface, and so $X$ admits a surjective morphism $\pi\colon
X\rightarrow\P^1$ whose generic fibre is irreducible.  The
N\'eron-Severi group has rank four, generated by the class $F$ of a
fibre of this morphism, the class $S$ of a section, and two more
classes, corresponding to irreducible components of reducible fibres,
or equivalently, exceptional divisors of the blowing down map to some
Hirzebruch surface.

There are three possible configurations of reducible fibres (see
Table~\ref{fibretable}):

\begin{table}
\begin{tabular}{ccc}
\scalebox{.5}{\includegraphics{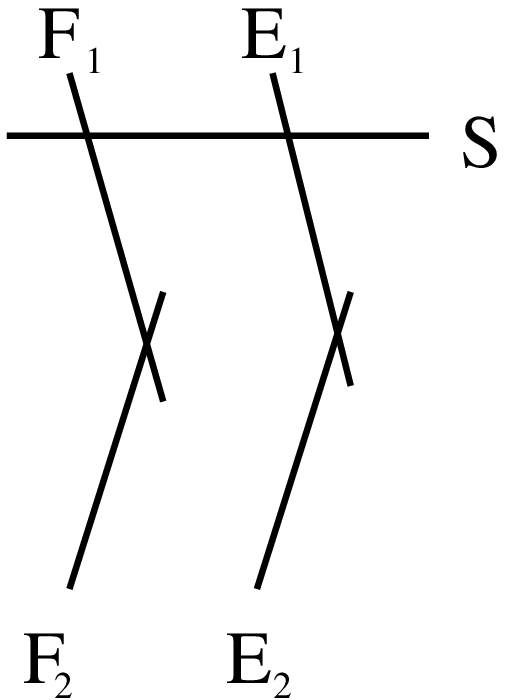}} &
\scalebox{.5}{\includegraphics{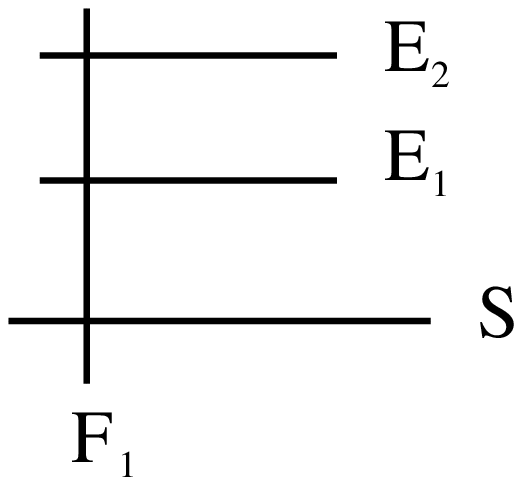}} &
\scalebox{.5}{\includegraphics{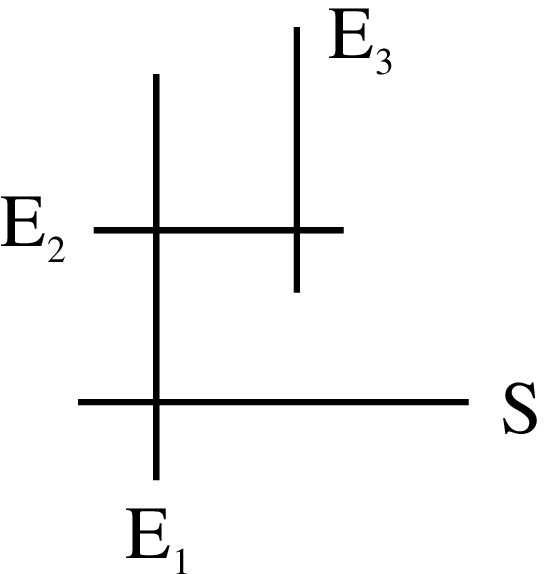}} \\
Case (1) & Case (2) & Case (3)
\end{tabular}
\vspace{.1in}
\caption{Configurations of reducible fibres}\label{fibretable}
\end{table}

\begin{enumerate}
\item Two reducible fibres, each with two components.  These components
intersect each other transversely in a single point.

\item One reducible fibre, with three components, configured like a
letter F, where the vertical component is the one which intersects the
section $S$.

\item One reducible fibre, with three components, configured like a 
letter H, where the leftmost vertical component is the one which intersects
the section $S$.

\end{enumerate}

In the proof, the first case will generate a further subcase,
corresponding to $\P^2$ blown up at three points in general position,
and the third case will generate a further subcase corresponding to a
double component in the reducible fibre.  However, in all cases, the
combinatorial description of the reducible fibres, combined with the
self-intersection of the section $S$, will completely determine the
intersection product on $X$.  Since our proof relies almost completely
on the intersection product, this will suffice to prove
Theorem~\ref{four}.

\vspace{.1in}

\noindent
Case (1): The fibration $\pi\colon X\rightarrow\P^1$ admits two reducible
fibres, each having two components.  (See Table~\ref{fibretable}.)

In addition to the classes $S$ and $F$ described above, let $E_1$ and
$E_2$ be the components of one irreducible fibre, and let $F_1$ and
$F_2$ be the components of the other.  Let $E_2$ and $F_2$ be the
components which do not intersect $S$.  Since the Picard rank is four,
the four classes $S$, $F$, $E_2$, and $F_2$ are a basis for the vector
space $NS(X)\otimes\R$.  The intersection matrix for $X$ is:
\[\begin{array}{r|rrrr}
& F & E_2 & F_2 & S \\ \hline
F & 0 & 0 & 0 & 1 \\
E_2 & 0 & -1 & 0 & 0 \\
F_2 & 0 & 0 & -1 & 0 \\
S & 1 & 0 & 0 & -n
\end{array}\]
where $n=-S^2$ is a positive integer.

The case $n>2$ is covered by Theorem~\ref{simplefibres}.  We first treat
the case $n=2$.

\begin{clm}
The effective cone of $X$ is generated by the classes of $E_i$, $F_i$,
and $S$, and the nef cone of $X$ is generated by the
classes of $F$, $D_2=2F+S$, $D_1=2F-E_2+S$, $D'_1=2F-F_2+S$, 
and $D_0=2F-E_2-F_2+S$.
\end{clm}

\noindent
{\it Proof:} \/ A straightforward calculation shows that the two cones
in the claim are dual to one another, and the first cone is clearly
contained in the effective cone, so it suffices to show that the interior
of the cone generated by $F$, $D_i$ and $D'_1$ is contained in the
ample cone.  We will use the Nakai-Moishezon Criterion again, and note 
that the intersection properties of the given five divisors are as follows:
\[\begin{array}{r|rrrrr}
& F & D_2 & D_1 & D'_1 & D_0 \\ \hline
F & 0 & 1 & 1 & 1 & 1 \\
D_2 & 1 & 2 & 2 & 2 & 2 \\
D_1 & 1 & 2 & 1 & 2 & 1 \\
D'_1 & 1 & 2 & 2 & 1 & 1 \\
D_0 & 1 & 2 & 1 & 1 & 0
\end{array}\]
Each of these five divisors are basepoint free: $F$ and $D_0$
correspond to morphisms to $\P^1$, $D_1$ and $D'_1$ correspond to
birational maps to $\P^2$, and $D_2=D_1+E_2=D_1'+F_2$.  It is
therefore clear that any positive linear combination of these five
divisors has positive self-intersection, and intersects any
irreducible curve $C$ positively, so by the Nakai-Moishezon Criterion,
must be ample.  The claim is therefore proven.  \qed

\vspace{.1in}

We will now prove Case (1) of Theorem~\ref{four} for the case $n=2$.
Let $P\in X(k)$ be any rational point --- we first assume that $P$
does not lie on $S$, $E_i$, or $F_i$.  Let $\psi\colon
X\rightarrow\P^1$ be the morphism corresponding to $D_0$.  A sequence
of best $D_0$-approximation to $P$ clearly lies along the fibre of
$\psi$ through $P$, while a sequence of best $D$-approximation to $P$
for any positive linear combination $D$ of the other four generators
lies along the component of $\pi$ through $P$.  We therefore define
two subcones of the nef cone of $X$: subcone $A$, generated by $F$,
$D_2$, $D_1$, $D'_1$, and $D_0+F$, and subcone $B$, generated by
$D_0$, $D_0+F$, $D_1$, and $D'_1$.  For each generator of subcone $A$,
a sequence of best approximation to $P$ lies along the component of
$F$ of minimal degree through $P$.  If $P$ does not lie on a reducible
fibre or $S$, then we can conclude from Corollary~\ref{divsum} that if
$D$ lies in the cone $A$, then the fibre of $\pi$ through $P$ is a
curve of best $D$-approximation to $P$.  For each generator of subcone
$B$, a similar argument shows that a sequence of best approximation to
$P$ lies along the component of the fibre of $\psi$ through $P$.
Since the union of $A$ and $B$ is the entire nef cone of $X$, we
conclude that Conjecture~\ref{ratcurve}is true for $P$.

If $P$ does lie on $S$, $E_i$, or $F_i$, then we divide the nef cone
into five subcones, one for each of the five curves $S$, $E_i$, and
$F_i$.  The subcone $A_C$ corresponding to a curve $C$ is the set of
nef divisors $D$ for which $C$ has minimal $D$-degree amongst $S$,
$E_i$, and $F_i$.  A straightforward calculation shows that generators
for these subcones can always be found from the set
$\{F,D_0,D_1,D'_1,D_2, F+D_0,F+D_1,F+D'_1,F+D_2+D_0\}$, and that each
divisor in that set intersects each of $S$, $E_i$, and $F_i$ in either
0 or 1.  It then follows from Corollary~\ref{divsum} that for each
curve $C\in\{S,E_1,E_2,F_1,F_2\}$, the curve $C$ is a curve of best
$D$-approximation to $P$ for every divisor $D\in A_C$.  Since the nef
cone is the union of the subcones $A_C$, we conclude that 
Conjecture~\ref{ratcurve} is true in case (1) with $n=2$.

We now move to the case $n=1$.  In this case, $X$ is the blowup of the
first Hirzebruch surface $H_1$ at two points, possibly ({\it a
priori}) infinitely near.  Thus, $X$ is also the blowup of $\P^2$ at
three points (again, possibly infinitely near).  If any of the three
points are infinitely near, then there will be a curve $S$ of
self-intersection $-2$ or less, and a morphism from $X$ to $\P^1$ of
which $S$ is a section, so $X$ is covered by the already-treated case
$n\geq 2$.  Thus, the three points must be distinct.  If they are
collinear, then no morphism from $X$ to $\P^1$ will have more than one
reducible section, so we cannot be in Case~(1).  Thus, $X$ must be
isomorphic to the blowup of $\P^2$ at three points in general
position.

Let $\phi\colon X\rightarrow\P^2$ be the blowing down map.  Let $L$ be
the divisor class corresponding to the invertible sheaf $\phi^*\O(1)$,
and let $E_i$ for $i=1,2,3$ be the three exceptional divisors of
$\phi$.  These four divisors are a basis of the vector space
$NS(X)\otimes\R$, and their intersection matrix is as follows:
\[\begin{array}{r|rrrr}
& L & E_1 & E_2 & E_3 \\ \hline
L & 1 & 0 & 0 & 0 \\
E_1 & 0 & -1 & 0 & 0 \\
E_2 & 0 & 0 & -1 & 0 \\
E_3 & 0 & 0 & 0 & -1 
\end{array}\]

\begin{clm}
The effective cone of $X$ is generated by the six curves $E_i$ and
$L-E_i-E_j$ for $i\neq j$.  The nef cone of $X$ is generated by the
five curves $L$, $L-E_i$, and $F=2L-E_1-E_2-E_3$.
\end{clm}

\noindent
{\it Proof:} \/ This can be found in, for example, \cite{Ma}.  \qed

\vspace{.1in}

The divisors $L$ and $F$ correspond to birational morphisms to $\P^2$,
and $L-E_i$ corresponds to a morphism to $\P^1$.  Thus, each of the
divisors is basepoint free, and their intersection matrix is:
\[\begin{array}{r|ccccc}
& L & L-E_1 & L-E_2 & L-E_3 & F \\ \hline
L & 1 & 1 & 1 & 1 & 2 \\
L-E_1 & 1 & 0 & 1 & 1 & 1 \\
L-E_2 & 1 & 1 & 0 & 1 & 1 \\
L-E_3 & 1 & 1 & 1 & 0 & 1 \\
F & 2 & 1 & 1 & 1 & 1
\end{array}\]

\vspace{.1in}

We now conclude the proof of Conjecture~\ref{ratcurve} in Case (1).
Let $P\in X(k)$ be any $k$-rational point, and let $D$ be any ample
divisor on $X$.  Let $A_1$ be the cone generated by the following
divisor classes:
\begin{equation}\label{conegens}
L,L-E_1,(L-E_1)+(L-E_2),(L-E_1)+(L-E_3),F
\end{equation}
By permuting the indices on the $E_i$ in the above list, we obtain two
more cones $A_2$ and $A_3$ such that the nef cone of $X$ is the union
of the subcones $A_i$.  Assume that $D$ lies in $A_1$' by symmetry, it
suffices to consider this case.

If $P$ does not lie on $E_i$ or $L-E_i-E_j$ for $i\neq j$, then
Theorem~\ref{projlines} and Corollary~\ref{divsum} imply that the
component of $L-E_m$ through $P$ is a curve of best $D$-approximation
to $P$.  (To see this, note that each of the generators of $A_1$ is
basepoint free, contracts no curves through $P$, and intersects
$L-E_1$ either one or zero times.)  Thus, if $P$ does not lie on any
$E_i$ or $L-E_i-E_j$, then we have proven Conjecture~\ref{ratcurve}
for $P$.

Now say that $P$ does lie on some $E_i$ or $L-E_i-E_j$.  The incidence
graph of these six divisors is a hexagon: $E_i$ intersects precisely
$L-E_i-E_k$ for $k\neq i$, and $L-E_i-E_j$ intersects precisely $E_i$
and $E_j$.  Thus, $P$ either lies on exactly one or exactly two of these
six curves.  

If $P$ lies only on $E_1$, then consider the subcone $B_1$ of $A_1$,
generated by:
\[L,L+(L-E_1),(L-E_1)+(L-E_2),(L-E_1)+(L-E_3),F\]
The generator $L$ contracts $E_1$, and so $E_1$ is clearly a curve of
best $L$-approximation to $P$.  The remaining generators are all
basepoint free, intersect $E_1$ once, and contract no curves through
$P$, and so by Theorem~\ref{projlines}, $E_1$ is a curve of best
approximation to $P$ with respect to all the generators of $B_1$.
Corollary~\ref{divsum} now implies that for any $D$ in $B_1$, $E_1$ is
a curve of best approximation to $P$.

If $D$ does not lie in $B_1$, then $D$ must lie in the cone $B'_1$,
generated by:
\[L+(L-E_1),L-E_1,(L-E_1)+(L-E_2),(L-E_1)+(L-E_3),F\]
The generator $L-E_1$ contracts the component $C$ of $L-E_1$ through
$P$, and so $C$ is a curve of best $(L-E_1)$-approximation to $P$.
Each of the other generators is basepoint free, does not contract any
curve through $P$, and intersects $C$ once, so by Theorem~\ref{projlines},
$C$ is a curve of best approximation to $P$ with respect to all the
generators of $B'_1$.  Corollary~\ref{divsum} now implies that $C$ is a
curve of best $D$-approximation to $P$.  In summary, if $P$ lies on $E_1$
and no other generator of the effective cone of $X$, then we have proven
Conjecture~\ref{ratcurve} for $P$.

If $P$ lies only on $E_2$, $E_3$, $L-E_1-E_2$, or $L-E_1-E_3$, then by
a similar argument to the previous paragraph, Theorem~\ref{projlines}
and Corollary~\ref{divsum} imply that for any $D$ in $A_1$, one of the
curves $E_2$, $E_3$, $L-E_1-E_2$, or $L-E_1-E_3$ will be a curve of
best $D$-approximation to $P$, and so Conjecture~\ref{ratcurve} is
proven for $P$.

If $P$ lies only on $L-E_2-E_3$, then define the cone $B_{23}$ to be 
that generated by:
\[F,L,(L-E_1)+(L-E_2),(L-E_1)+(L-E_3),F+(L-E_1)\]
Since $F$ contracts $L-E_2-E_3$, it's clear that $L-E_2-E_3$ is a
curve of best $F$-approximation to $P$.  All the other generators are
basepoint free, do not contract any curves through $P$, and intersect
$L-E_2-E_3$ once.  Thus, by Theorem~\ref{projlines} and
Corollary~\ref{divsum}, we conclude that if $D$ lies in $B_{23}$,
$L-E_2-E_3$ is a curve of best $D$-approximation to $P$.  If $D$ does
not lie in $B_{23}$, then $D$ must lie in $B'_{23}$, generated by:
\[L-E_1,L,(L-E_1)+(L-E_2),(L-E_1)+(L-E_3),F+(L-E_1)\]
A similar argument to the previous paragraph shows that the component
of $L-E_1$ through $P$ is a curve of best $D$-approximation to $P$.  We
conclude that if $P$ lies only on one of the six generators of the effective
cone of $X$, then Conjecture~\ref{ratcurve} is true for $P$.

It only remains to consider the possibility that $P$ is the intersection
of two of the curves $E_i$ or $L-E_i-E_j$.  First, assume that $P=E_1\cap
L-E_1-E_i$ for $i=2$ or $i=3$.  Let $B_{11i}$ be the cone generated by:
\[L+(L-E_1),L+F,L+(L-E_1)+(L-E_i),(L-E_1)+(L-E_{5-i}),\]
\[L-E_1,(L-E_1)+(L-E_i),F\]
The last three generators of $B_{11i}$ all contract $L-E_1-E_i$, and so
$L-E_1-E_i$ is a curve of best approximation to $P$ for each of these
generators.  The first four generators are all basepoint free,
contract no curves through $P$, and map $L-E_1-E_i$ to a line, and so
by Theorem~\ref{projlines} and Corollary~\ref{divsum} we conclude that
$L-E_1-E_2$ is a curve of best $D$-approximation to $P$ if $D$ lies in
$B_{11i}$.

If $D$ does not lie in $B_{11i}$, then $D$ must lie in the cone $B'_{11i}$,
generated by:
\[L,L+(L-E_1),L+(L-E_1)+(L-E_i),(L-E_1)+(L-E_{5-i}),L+F\]
(Recall that $D$ is assumed to lie in the cone $A_1$, defined above.)
The divisor $L$ contracts $E_1$, and the other four generators are
basepoint free, contract no curves through $P$, and map $E_1$ to a
line.  We conclude from Theorem~\ref{projlines} and
Corollary~\ref{divsum} that $E_1$ is a curve of best $D$-approximation
to $P$.

Next, assume that $P$ lies on the intersection of $L-E_1-E_i$ and $E_i$
for $i=2$ or $i=3$.  Define the cone $B_{1ii}$ to be the cone generated by:
\[L,L-E_1,(L-E_1)+(L-E_{5-i}),L+(L-E_1)+(L-E_i),L+F,\]
\[F+(L-E_1)+(L-E_{5-i})\]
All of these classes are basepoint free.  The first three contract
$E_i$, and the last three contract no curves through $P$ and map $E_i$
to a line, and so we conclude by Theorem~\ref{projlines} and
Corollary~\ref{divsum} that $E_i$ is a curve of best $D$-approximation
to $P$ for any $D$ in $B_{1ii}$.  Thus, Conjecture~\ref{ratcurve} is
true for $P$ if $D$ lies in $B_{1ii}$.

If $D$ does not lie in $B_{1ii}$, then $D$ must lie in the cone $B'_{1ii}$,
generated by:
\[F,L-E_1,(L-E_1)+(L-E_i),L+F,L+(L-E_1)+(L-E_i),\]
\[F+(L-E_1)+(L-E_{5-i})\]
All generators are basepoint free.  The first three generators
contract $L-E_1-E_i$, and the last three contract no curves through
$P$ but map $L-E_1-E_i$ to a line.  We conclude by
Theorem~\ref{projlines} and Corollary~\ref{divsum} that $L-E_1-E_i$
is a curve of best $D$-approximation to $P$ for all $D$ in $B'_{1ii}$.
We conclude that Conjecture~\ref{ratcurve} is true for
$P=E_i\cap(L-E_1-E_i)$ and all $D$.

Finally, let $P$ be the intersection of $E_i$ and $L-E_i-E_{5-i}$ for 
$i=2$ or $i=3$.  Define $B_{123}$ to be the cone generated by:
\[L,L-E_1,(L-E_1)+(L-E_{5-i}),(L-E_1)+(L-E_i),F+L,F+(L-E_1),\]
\[F+(L-E_1)+(L-E_{5-i})\]
All of these classes are basepoint free.  The first three contract $E_i$,
and the last four do not contract any curve through $P$, and map $E_i$ to
a line.  Thus, by Theorem~\ref{projlines} and Corollary~\ref{divsum},
we conclude that $E_i$ is a curve of best $D$-approximation to $P$, 
provided $D$ lies in $B_{123}$.

If $D$ does not lie in $B_{123}$, then it must lie in the cone $B'_{123}$,
generated by:
\[F,F+L,(L-E_1)+(L-E_i),F+(L-E_1),F+(L-E_1)+(L-E_{5-i})\]
All five of these are basepoint free.  The last four contract no
curves through $P$, and map $L-E_2-E_3$ to a line, and $F$ contracts
$L-E_2-E_3$ to a point.  By Theorem~\ref{projlines} and
Corollary~\ref{divsum}, we conclude that $L-E_i-E_j$ is a curve of
best $D$-approximation to $P$ if $D$ lies in $B'_{123}$.
Conjecture~\ref{ratcurve} is thus proven for $X$.

\vspace{.1in}

Case (2): One reducible fibre, with three components, configured like a
letter F, where the vertical component is the one which intersects the
section $S$.  (See Table~\ref{fibretable}.)

\vspace{.1in}

Let $F_1$ be the component of the reducible fibre which intersects $S$,
and let $E_1$ and $E_2$ be the other two components.  Let $F$ be the class
of a fibre.  These four classes are a basis of $\mbox{NS}(X)$, with
intersection matrix:
\[\begin{array}{r|rrrr}
& S & E_1 & E_2 & F \\ \hline
S & -n & 0 & 0 & 1 \\
E_1 & 0 & -1 & 0 & 0 \\
E_2 & 0 & 0 & -1 & 0 \\
F & 1 & 0 & 0 & 0 
\end{array}\]

\begin{clm}
The curves $S$, $F_1$, $E_1$, and $E_2$ generate the effective cone of
$X$.  The nef cone of $X$ is generated by the divisors $F$, 
$D_1=nF+S$, $D_2=nF+S-E_1$, and $D_3=nF+S-E_2$.
\end{clm}

\noindent
{\it Proof:} \/ A straightforward calculation shows that the two cones
described in the claim are dual to one another.  Thus, it suffices to
show that $F$, $D_1$, $D_2$, and $D_3$ generate the nef cone of $X$.
To do this, it further suffices to show that every positive linear
combination of $F$, $D_1$, $D_2$, and $D_3$ is ample, since the converse
inclusion is clear from the effectivity of $S$, $F_1$, $E_1$, and $E_2$.

The intersections of $F$ and the $D_i$ are tabulated as follows:
\[\begin{array}{r|cccc}
& F & D_1 & D_2 & D_3 \\ \hline
F & 0 & 1 & 1 & 1 \\
D_1 & 1 & n & n & n \\
D_2 & 1 & n & n-1 & n \\
D_3 & 1 & n & n & n-1 
\end{array}\]

Since $n\geq 1$, it is clear that any positive linear combination of
$F$ and the $D_i$ must have positive self-intersection.  It's clear
that $F$ is basepoint free, since it corresponds to a morphism to
$\P^1$, and similarly each $D_i$ corresponds to a morphism to a cone
over a rational normal curve (see Remark~\ref{actualcones}).  In
particular, each of these four divisors is nef, so by the
Nakai-Moishezon criterion for ampleness, we conclude that they must
generate the nef cone of $X$, as desired.  \qed

\vspace{.1in}

\begin{rem}\label{actualcones}
Let $Y$ be the cone over a smooth rational normal curve of degree $n$.
It is easy to see that if we blow up $Y$ at the vertex, we obtain the
Hirzebruch surface $H_n$.  Moreover, the exceptional curve of the
blowup $\pi\colon H_n\to Y$ is precisely the unique $(-n)$-section $S'$
of $H_n$, and if $L$ is the class of a hyperplane section of $Y$, then
$\pi^*L=S'+nF'$, where $F'$ is a fibre of the morphism $H_n\to\P^1$.

Our rational surface $X$ admits several birational maps to Hirzebruch
surfaces.  It is a straightforward calculation to see that each $D_i$
is the pullback of $\pi^*L$ via one of these maps.  For example, we
might blow down $E_1$ and then $F_1$ to obtain $H_n$, in which case
pulling back $\pi^*L$ to $X$ gives $D_3$.  This technique will be used
frequently in what follows.
\end{rem}

\vspace{.1in}

We now prove Conjecture~\ref{ratcurve} for Case (2).  We first assume
that $n>1$.  Let $P\in X(k)$ be any rational point not lying on $S$,
$F_1$, or either $E_i$.  Let $C$ be the (irreducible) component of $F$
through $P$.  Then $F$ contracts $C$, and all of the other generators
of the nef cone map $C$ to a line, and so by Theorem~\ref{projlines}
and Corollary~\ref{divsum}, we conclude that $C$ is a curve of best
$D$-approximation to $P$ for any ample divisor $D$.

If $P$ lies on $S$ but not $F_1$, then we divide the nef cone of $X$
into two subcones: a cone $A$ generated by $D_1$, $D_2$, $D_3$,
$D_1+F$, $D_2+F$, and $D_3+F$, and a cone $B$ generated by $D_1+F$,
$D_2+F$, $D_3+F$, and $F$.  The classes $D_i$ each contract $S$, and
$F+D_i$ is basepoint free, contracts no curves through $P$ but maps
$S$ to a line, and so for each of these six divisor classes, $S$ is a
curve of best approximation to $P$.  Thus, by Corollary~\ref{divsum},
Conjecture~\ref{ratcurve} is true for $P$ with respect to any ample
divisor in $A$.  On the other hand, $F+D_i$ also maps the component
$C$ of $F$ through $P$ to a line, and $F$ contracts $C$, so similar
reasoning shows that Conjecture~\ref{ratcurve} is true for $P$ with
respect to any ample divisor in $B$.  We conclude that
Conjecture~\ref{ratcurve} is true for $P$.

If $P$ lies on $E_i$ but not $F_1$, or if $P$ lies on $F_1$ but not
$S$ or either $E_i$, then let $C$ be the unique irreducible component
of $F$ through $P$ (that is, either $F_1$ or $E_i$).  Every generator
of the nef cone either contracts $C$ or else maps $C$ to a line
and contracts no curves through $P$, and so Theorem~\ref{projlines}
and Corollary~\ref{divsum} imply as usual that $C$ is a curve of
best $D$-approximation to $P$ for any ample divisor $D$.

If $P=S\cap F_1$, then consider the cone $A'$ generated by:
\[F,F+D_1,D_2,D_3\]
Note that all four of these classes are basepoint free.  All but
$F+D_1$ contract $F_1$, and $F+D_1$ maps $F_1$ to a line and contracts
no curve through $P$.  We conclude from Theorem~\ref{projlines} and
Corollary~\ref{divsum} that $F_1$ is a curve of best $D$-approximation 
to $P$ for any $D$ in $A'$.  If $D$ is not in $A'$, then $D$ must lie in 
$B'$, generated by:
\[D_1,F+D_1,D_2,D_3\]
Every generator is basepoint free, and all but $F+D_1$ contracts $S$
The class $F+D_1$ maps $S$ to a line, and contracts no curves through $P$.
We conclude by Theorem~\ref{projlines} and Corollary~\ref{divsum} that
$S$ is a curve of best $D$-approximation to $P$ for every $D$ in $B'$.
Thus, Conjecture~\ref{ratcurve} is true for $P$.

Finally, if $P=F_1\cap E_i$ for some $i$, then consider the cone
$\tilde{A}$ generated by $F$, $D_1+D_{i+1}$, $D_2$, and $D_3$.  Each
generator is basepoint free, and all but $D_1+D_{i+1}$ contract $F_1$.
Furthermore, the class $D_1+D_{i+1}$ maps $F_1$ to a line and
contracts no curves through $P$.  Theorem~\ref{projlines} and
Corollary~\ref{divsum} now imply that $F_1$ is a curve of best
$D$-approximation to $P$ for any $D$ in $\tilde{A}$.

If $D$ does not lie in $\tilde{A}$, then it must lie in $\tilde{B}$,
generated by $F$, $D_1$, $D_1+D_{i+1}$, and $D_{4-i}$.  Each generator
is basepoint free, and all but $D_1+D_{i+1}$ contract $E_i$.
Furthermore, the class $D_1+D_{i+1}$ maps $E_i$ to a line and
contracts no curves through $P$.  Theorem~\ref{projlines} and
Corollary~\ref{divsum} now imply that $E_i$ is a curve of best
$D$-approximation to $P$ for any $D$ in $\tilde{B}$.  We conclude that
Conjecture~\ref{ratcurve} is true for $X$ in the case that $n>1$.

(Note that if $n>3$, this result also follows immediately from
Theorem~\ref{onefibre}.)

If $n=1$, then $D_2^2=D_3^2=0$, so these both correspond to morphisms
to $\P^1$.  In fact, if we cannot describe $X$ with $n>1$, then $X$
must be a blowup of $\P^1\times\P^1$ at two (possibly infinitely near)
points.  In fact, the points cannot be infinitely near, since otherwise
the exceptional divisor of the first blowup would have self-intersection
$-2$, and we could choose it to be the section $S$, with $n=2>1$.  And
if the configuration of exceptional divisors is to be as in Case (2),
the two blown up points must lie on the same fibre of at least one
of the canonical projections $\P^1\times\P^1\rightarrow\P^1$.

In other words, we have shown that $X$ must be isomorphic to $\P^2$
blown up at three different points on a straight line.  The
exceptional curves of this blowup are $S$, $E_1$, and $E_2$, while
$F_1$ is the strict transform of the line joining the blown up points.
The divisors $D_2$ and $D_3$ correspond to the maps $\pi_i\colon
X\rightarrow\P^1$ ($i=2,3$) which factor through the blowup of $\P^2$
at a single point, while $F$ corresponds to the third map $\pi_1\colon
X\rightarrow\P^1$, induced in the same way.  The divisor $D_1$
corresponds to the map $\phi\colon X\rightarrow\P^2$ that blows down
$S$, $E_1$, and $E_2$.

Let $P\in X(k)$ be a rational point which does not lie on $S$, $F_1$, or
either $E_i$, and let $D$ be an ample divisor on $X$.  By symmetry, we
may assume that $D.F\leq D.D_i$ for $i=2,3$, since $F$, $D_2$, and $D_3$
are all conjugate under the automorphism group of $X$.  In that case, $D$
lies in the cone $A$ generated by:
\[F,D_1,F+D_2,F+D_3,F+D_2+D_3\]
All these divisors are basepoint free.  The last four contract no
curves through $P$ and map the component $C$ of $F$ through $P$ to a
line, and $F$ contracts $C$.  Thus, by Theorem~\ref{projlines} and
Corollary~\ref{divsum}, we conclude that $C$ is a curve of best
$D$-approximation to $P$, and that Conjecture~\ref{ratcurve} is true
for $P$.

If $P$ lies on one or more of $S$, $E_1$, $E_2$, or $F_1$, then we must
subdivide $A$ further.  If $P$ lies on $S$ but not $F_1$, then assume
that $D$ lies in the cone $A_S$ generated by:
\[D_1,F+D_1,F+D_2,F+D_3,F+D_2+D_3\]
Each of these divisors is basepoint free.  The last four contract no
curves through $P$, and map $S$ to a line, and $D_1$ contracts $S$, so
by Theorem~\ref{projlines} and Corollary~\ref{divsum}, we conclude that 
$S$ is a curve of best $D$-approximation to $P$.  If $D$ does not lie in
$A_S$, then $D$ must lie in $A'_S$, generated by:
\[F,F+D_1,F+D_2,F+D_3\]
All of these divisors are basepoint free.  The last four contract no
curves through $P$ and map the component $C$ of $F$ through $P$ to a
line, and $D_1$ contracts $C$.  We conclude from
Theorem~\ref{projlines} and Corollary~\ref{divsum} that $C$ is a curve
of best $D$-approximation to $P$.  In all cases, if $P$ lies on $S$
but not $F_1$, then Conjecture~\ref{ratcurve} is true for $P$.

If $P$ lies on exactly one of $E_1$, $E_2$, or $F_1$, then a quick
check (using Theorem~\ref{projlines} and Corollary~\ref{divsum}) shows
that $E_i$ is a curve of best $D$-approximation for every $D$ in the
cone $A$.  The same is true if $P=S\cap F_1$.

The only remaining case is thus $P=F_1\cap E_i$.  Consider the cone
$A_i$ generated by:
\[F,D_1,F+D_{4-i},F+D_1+D_{i+1},F+D_1+D_2+D_3\]
All these divisors are basepoint free.  The first three generators
contract $E_i$, while the last two contract no curves through $P$ but
map $E_i$ to a line.  By Theorem~\ref{projlines} and Corollary~\ref{divsum},
we conclude that $E_i$ is a curve of best $D$-approximation to $P$ for
any $D$ in $A_i$.  If $D$ does not lie in $A_i$, then $D$ must lie in the
cone $A'_i$, generated by:
\[F,F+D_2,F+D_3,F+D_2+D_3,F+D_1+D_2,F+D_1+D_2+D_3\]
All of these divisors are basepoint free.  The first four contract $F_1$,
and the last two contract no curves through $P$ but map $F_1$ to a line,
so we conclude by Theorem~\ref{projlines} and Corollary~\ref{divsum} that
$F_1$ is a curve of best $D$-approximation to $P$.

The sum of all these calculations is that Conjecture~\ref{ratcurve} is
true for $X$ in Case~(2).

\vspace{.1in}

Case (3): One reducible fibre, with three components, configured like
a letter H, where the leftmost vertical component is the one which
intersects the section $S$.  (See Table~\ref{fibretable}.)

\vspace{.1in}

We first deal with the case in which the reducible fibre has no
multiple components.  Let $E_1$ be the component of the reducible
fibre which intersects the section $S$.  Let $E_2$ be the component
which intersects both other components, and let $E_3$ be the remaining
component.  Let $F=E_1+E_2+E_3$ be the class of a fibre.  These four
classes are a basis of $\mbox{NS}(X)$, with intersection matrix:
\[\begin{array}{r|rrrr}
& S & E_2 & E_3 & F \\ \hline
S & -n & 0 & 0 & 1 \\
E_2 & 0 & -2 & 1 & 0 \\
E_3 & 0 & 1 & -1 & 0 \\
F & 1 & 0 & 0 & 0 
\end{array}\]

\begin{clm}
The curves $S$, $E_1$, $E_2$, and $E_3$ generate the effective cone of
$X$.  The nef cone of $X$ is generated by the divisors $F$, 
$D_1=nF+S-E_2-2E_3$, $D_2=nF+S-E_2-E_3$, and $D_3=nF+S$.
\end{clm}

\noindent
{\it Proof:} \/ A straightforward calculation shows that the two cones
described in the claim are dual to one another.  Thus, it suffices to
show that $F$, $D_1$, $D_2$, and $D_3$ generate the nef cone of $X$.
To do this, it further suffices to show that every positive linear
combination of $F$, $D_1$, $D_2$, and $D_3$ is ample, since the converse
inclusion is clear from the effectivity of $S$, $F_1$, $E_1$, and $E_2$.

The intersections of $F$ and the $D_i$ are tabulated as follows:
\[\begin{array}{r|cccc}
& F & D_1 & D_2 & D_3 \\ \hline
F & 0 & 1 & 1 & 1 \\
D_1 & 1 & n-2 & n-1 & n \\
D_2 & 1 & n-1 & n-1 & n \\
D_3 & 1 & n & n & n 
\end{array}\]

We first note that $n$ can be chosen to be at least two, since if
$n=0$ or $n=1$, then $X$ is a blowup of $\P^2$ at two infinitely near
points, and is therefore also a blowup of the second Hirzebruch
surface $H_2$, in which case we can take $n=2$.

Now assume that $n\geq 2$.  It's clear that $F$ and each of the $D_i$
are basepoint free: $F$ corresponds to a morphism to $\P^1$, $D_1$
corresponds to a morphism to a cone over a smooth rational curve of
degree $n-2$ (if $n=2$ we take this cone to be $\P^1$), $D_2$
corresponds to a morphism to a cone over a smooth rational curve of
degree $n-1$, and $D_3$ corresponds to a morphism to a cone over a
smooth rational curve of degree $n$ (see Remark~\ref{actualcones}).
For each $D_i$, the ruling of the cone pulls back to $F$.  In
particular, every positive linear combination of $F$ and the $D_i$ has
positive intersection with every curve on $X$, and has positive
self-intersection.  Therefore, by the Nakai-Moishezon criterion, every
such positive linear combination is ample, and we have proven the
claim.  \qed

\vspace{.1in}

We now prove Conjecture~\ref{ratcurve} for $X$.  Let $P\in X(k)$ be
any rational point.  The case $n>3$ is immediate from
Theorem~\ref{onefibre}.

Assume $n=3$, and let $P\in X(k)$ be any $k$-rational point.  If $P$
does not lie on $S$ or any $E_i$, then every generator of the nef cone
either contracts the component $C$ of $F$ through $P$, or else
contracts no curves through $P$ and maps $C$ to a line.  Thus,
Theorem~\ref{projlines} and Corollary~\ref{divsum} imply that for any
ample divisor $D$, $C$ is a curve of best $D$-approximation to $P$.

If $P$ lies on $S$ but not $E_1$, then consider the cone $A_S$, generated
by:
\[F,F+D_1,F+D_2,F+D_3\]
Each divisor is basepoint free.  The last three divisors contract no
curves through $P$, but map the component $C$ of $F$ through $P$ to a
line, while $F$ contracts $C$.  By Theorem~\ref{projlines} and
Corollary~\ref{divsum}, we conclude that $C$ is a curve of best
$D$-approximation to $P$ for any $D$ in $A_S$.  If $D$ does not lie in
$A_S$, then it must lie in $A'_S$, generated by:
\[D_1,F+D_1,D_2,F+D_2,D_3,F+D_3\]
Each $D_i$ contracts $S$, while $F+D_i$ contracts no curves through
$P$ and maps $S$ to a line.  Thus, by Theorem~\ref{projlines} and
Corollary~\ref{divsum}, we conclude that $S$ is a curve of best
$D$-approximation to $P$ for every $D\in A'_S$.  Conjecture~\ref{ratcurve}
is thus proven for $P$.

If $P$ lies on exactly one of the $E_i$ but not $S$, then for each
generator $G$ of the nef cone, we see that either $G$ contracts $E_i$
or else maps $E_i$ to a line and contracts no curves through $P$.
Theorem~\ref{projlines} and Corollary~\ref{divsum} now imply that
$E_i$ is a curve of best $D$-approximation to $P$.

If $P=S\cap E_1$, then consider the cone $A_1$, generated by:
\[F,D_1,D_2,F+D_3\]
The first three all contract $E_1$, and the last contracts no curves
through $P$ but maps $E_1$ to a line.  By Theorem~\ref{projlines} and
Corollary~\ref{divsum}, we conclude that $E_1$ is a curve of best
$D$-approximation to $P$ for all $P\in A_1$.  If $D$ does not lie in $A_1$,
then $D$ must lie in $A'_1$, generated by:
\[D_1,D_2,D_3,F+D_3\]
The first three curves all contract $S$, while the last contracts no
curves through $P$ but maps $S$ to a line.  Theorem~\ref{projlines}
and Corollary~\ref{divsum} now imply that $S$ is a curve of best
$D$-approximation to $P$ for any $D$ in $A'_1$.

If $P=E_1\cap E_2$, then consider the cone $A_2$, generated by:
\[F,D_1,D_3,D_2+D_3\]
The first three classes all contract $E_2$, and the last contracts no
curves through $P$ but maps $E_2$ to a line.  We conclude by
Theorem~\ref{projlines} and Corollary~\ref{divsum} that $E_2$ is a
curve of best $D$-approximation to $P$ for all $D$ in $A_2$.  If $D$
does not lie in $A_2$, then $D$ must lie in $A'_2$, generated by:
\[F,D_1,D_2,D_2+D_3\]
The first three contract $E_1$, and the last contracts no curve
through $P$ but maps $E_1$ to a line.  By Theorem~\ref{projlines} and
Corollary~\ref{divsum}, we conclude that $E_1$ is a curve of best
$D$-approximation to $P$ for all $D$ in $A'_2$.
Conjecture~\ref{ratcurve} is thus proven for $P$.

Finally, if $P=E_2\cap E_3$, then consider the cone $A_3$, generated by:
\[F,D_2,D_3,D_1+D_2\]
The first three all contract $E_3$, while the last contracts no curves
through $P$ but maps $E_3$ to a line.  We conclude by
Theorem~\ref{projlines} and Corollary~\ref{divsum} that $E_3$ is a
curve of best $D$-approximation to $P$ for any $D$ in $A_3$.  If $D$
does not lie in $A_3$, then $D$ must lie in $A'_3$, generated by:
\[F,D_1,D_3,D_2+D_3\]
The first three contract $E_2$, while that last contracts no curves
through $P$ but maps $E_2$ to a line.  Theorem~\ref{projlines} and
Corollary~\ref{divsum} now imply that $E_2$ is a curve of best
$D$-approximation to $P$ for all $D$ in $A'_3$.  We conclude that
Conjecture~\ref{ratcurve} is true for this $P$.  Since this exhausts
the possibilities for $P$, we conclude that Conjecture~\ref{ratcurve}
is proven for $X$ in Case~(3) with $n=3$.

If $n=2$, then both $F$ and $D_1$ correspond to maps to $\P^1$, so
that $X$ is realized as the blowup of $\P^1\times\P^1$ at two
infinitely near points.  Let $P$ be any $k$-rational point on $X$
which does not lie on $S$ or any $E_i$.  Let $C_1$ and $C_F$ be the
components of $D_1$ and $F$, respectively, through $P$.  Consider the
cone $A$, generated by:
\[F,F+D_1,D_2,D_3,D_1+D_3\]
All of these classes are basepoint free, and all but the last either
contract $C_F$ or else contract no curves through $P$ and map $C_F$ to
a line.  Thus, if we can show that $C_F$ is a curve of best
$(D_1+D_3)$-approximation to $P$, then Theorem~\ref{projlines} and
Corollary~\ref{divsum} will imply that $C_F$ is a curve of best
$D$-approximation to $P$ for any $D$ in $A$.

To prove this, note that $C_F.D_1=C_F.D_3=1$, so that a sequence of
best $D_1$- or $D_3$-approximation to $P$ along $C.F$ has constant of
approximation equal to 1.  Thus, by Theorem~\ref{product}, since $D_1$
and $D_3$ are both basepoint free, it follows that any sequence with a
better constant of $(D_1+D_3)$-approximation to $P$ than $1+1=2$ must
have all but finitely many points contained in a union of curves
through $P$ which are contracted by $D_1$ or $D_3$.  The only such
curve is $C_1$, which is only contracted by $D_1$.  Since
$C_1.(D_1+D_3)=2$, it follows that no sequence of points on $C_1$ can
have a constant of $(D_1+D_3)$-approximation to $P$ of better than 2.
Thus, $C_F$ is a curve of best $(D_1+D_3)$-approximation to $P$, and
since $C_1.(D_1+D_3)=2$, we see that $C_1$ is also a curve of best
$(D_1+D_3)$-approximation to $P$.  (Recall that $P$ is assumed not to
lie on $S$ or any $E_i$.)  We conclude that $C_F$ is a curve of best
$D$-approximation to $P$ for any $D$ in $A$.

If $D$ does not lie in $A$, then it must lie in the cone $B$, generated by:
\[D_1,F+D_1,D_2,D_1+D_3\]
Each divisor is basepoint free, and the first three each either
contract $C_1$ or else contract no curves through $P$ and map $C_1$ to
a line, and we have already seen that $C_1$ is a curve of best
$(D_1+D_3)$-approximation to $P$.  Thus, by Corollary~\ref{divsum}, we
conclude that $C_1$ is a curve of best $D$-approximation to $P$ for
any $D$ in $B$.  It follows that Conjecture~\ref{ratcurve} is true for
$P$.

If $P$ lies on $S$ but not any $E_i$, then for any $D$ in the cone $B$,
it is easy to check that $S$ is a curve of best $D$-approximation to $P$.
If $D$ lies in $A$, then consider the cone $A_S$, generated by:
\[D_1+D_3,D_2,D_3,F+D_1,F+D_2,F+D_3\]
The first three of these contract $S$, and the last three contract no
curves through $P$ but map $S$ to a line.  By Theorem~\ref{projlines}
and Corollary~\ref{divsum}, it follows that $S$ is a curve of best
$D$-approximation to $P$ for all $D$ in $A_S$.  If $D$ does not lie in 
$A_S$ or $B$, then it must lie in the cone $A'_S$, generated by:
\[F,F+D_1,F+D_2,F+D_3\]
The first of these contracts $C_F$, and the last three contract no
curves through $P$ but map $C_F$ to a line.  By
Theorem~\ref{projlines} and Corollary~\ref{divsum}, it follows that
$C_F$ is a curve of best $D$-approximation to $P$ for all $D$ in
$A'_S$.  We conclude that Conjecture~\ref{ratcurve} is true for $P$.

Next, assume that $P$ lies on a curve $C$ which is either $E_1$ or $E_2$,
and assume that $P$ does not lie on $S$ or $E_3$, and is not the point 
$E_1\cap E_2$.  It is easy to check that each generator of the nef cone of
$X$ either contracts $C$ or else contracts no curves through $P$ and maps
$C$ to a line.  We conclude that for any ample divisor $D$, the curve $C$
is a curve of best $D$-approximation to $P$.

If $P$ lies on $E_3$ but not $E_2$, then it is easy to check that
every generator of the cone $A$ either contracts $E_3$ or else
contracts no curves through $P$ but maps $E_3$ to a line.  Thus, if
$D$ lies in $A$, then $E_3$ is a curve of best $D$-approximation to
$P$.  If $D$ does not lie in $A$, then consider the cone $A_3$,
generated by:
\[D_2,D_1+D_3,F+D_1,D_1+D_2,2D_1+D_3\]
It is clear that $E_3$ is contracted by $D_2$ and mapped to a line by
the next three divisors (which contract no curves through $P$).
Moreover, $E_3$ is a curve of best $D_1$-approximation to $P$, and
therefore is also a curve of best $2D_1$-approximation to $P$, and by
Theorem~\ref{addeffective}, since $E_3.D_3=0$, it follows that $E_3$
is also a curve of best $(2D_1+D_3)$-approximation to $P$.  By
Corollary~\ref{divsum}, it follows that $E_3$ is a curve of best
$D$-approximation to $P$ for any $D$ in $A_3$.  If $D$ does not lie in
either $A_3$ or $A$, then it must lie in $B_3$, generated by:
\[D_1,F+D_1,D_1+D_2,2D_1+D_3\]
The first of these contracts $C_1$, the next two contract no curves
through $P$ but map $C_1$ to a line, and since
$(2D_1+D_3).C_1=2=(2D_1+D_3).E_3$, it follows that $C_1$ is a curve of
best $(2D_1+D_3)$-approximation to $P$ as well.  Thus, by
Corollary~\ref{divsum}, we deduce that $C_1$ is a curve of best
$D$-approximation to $P$ for any $D$ in $B_3$, and hence for any $D$
in the nef cone.

This leaves just three points $P$ to check.  If $P=S\cap E_1$, then 
consider the cone $A_1$, generated by:
\[F+D_3,D_1,D_2,D_3\]
All of these are basepoint free, and all but the first of them contract
$S$.  Since $F+D_3$ contracts no curves through $P$ but maps $S$ to a 
line, we conclude by Theorem~\ref{projlines} and Corollary~\ref{divsum}
that $S$ is a curve of best $D$-approximation to $P$ for all $D$ in $A_1$.
If $D$ does not lie in $A_1$, then $D$ must lie in $B_1$, generated by:
\[F,D_1,D_2,F+D_3\]
The first three of these contract $E_1$, and the last maps $E_1$ to a
line but does not contract any curve through $P$.  By
Theorem~\ref{projlines} and Corollary~\ref{divsum}, we conclude that
$E_1$ is a curve of best $D$-approximation to $P$ for any $D$ in
$B_1$, and hence that Conjecture~\ref{ratcurve} is true for $P$.

If $P=E_1\cap E_2$, then consider the cone $A_2$, generated by:
\[F,D_1,D_2,D_2+D_3\]
The first three of these contract $E_1$, and the last contracts no curves
through $P$ but maps $E_1$ to a line.  We conclude by Theorem~\ref{projlines}
and Corollary~\ref{divsum} that $E_1$ is a curve of best $D$-approximation
to $P$ for all $D$ in $A_2$.  If $D$ does not lie in $A_2$, then $D$ must lie
in $B_2$, generated by:
\[F,D_1,D_3,D_2+D_3\]
By an exactly similar argument to that used for $A_2$, we conclude
that $E_2$ is a curve of best $D$-approximation to $P$ for any $D$ in
$B_2$, and therefore that Conjecture~\ref{ratcurve} is true for $P$.

Finally, we consider the case that $P=E_2\cap E_3$.  Let $A'_2$ be
the cone generated by:
\[F,D_1,D_3,D_1+D_2\]
The first three divisors contract $E_2$, and the last contracts no
curves through $P$ but maps $E_2$ to a line.  Theorem~\ref{projlines}
and Corollary~\ref{divsum} thus imply that $E_2$ is a curve of best
$D$-approximation to $P$ for any $D$ in $A'_2$.  If $D$ does not lie in
$A'_2$, then it must lie in $B'_2$, generated by:
\[F,D_2,D_3,D_1+D_2\]
Again, each of the first three divisors contracts $E_3$, while the
last contracts no curves through $P$ but maps $E_3$ to a line.  We
conclude from Theorem~\ref{projlines} and Corollary~\ref{divsum} that
$E_3$ is a curve of best $D$-approximation to $P$ for any $D$ in
$B'_2$, and therefore that Conjecture~\ref{ratcurve} is true for $P$.
This completes the proof of Conjecture~\ref{ratcurve} for Case~(3),
provided that the reducible fibre has no multiple components.

To finish the proof, it remains only to consider the case of a
multiple component in the reducible fibre.  The only way this can
occur is if the crossbar $E_2$ of the H has multiplicity two in the
fibre.

\begin{clm}
The effective cone of $X$ is generated by the classes of $S$, $E_1$,
$E_2$, and $E_3$.  The nef cone of $X$ is generated by the classes
of $F$, $D_1=S+nF$, $D_2=2S+2nF-2E_2-E_3$, and $D_3=S+nF-E_2-E_3$.
\end{clm}

\noindent
{\it Proof of claim:} \/ A simple calculation shows that the two cones
are dual, so to prove the claim it suffices to show that every
positive linear combination of $F$, $D_1$, $D_2$, and $D_3$ is ample.
The classes $D_1$ and $D_3$ are basepoint free, since they correspond
to morphisms to cones (see Remark~\ref{actualcones}), and $F$ is also
basepoint free because it corresponds to a morphism to $\P^1$.
Finally, note that $D_2$ is basepoint free because it is linearly
equivalent both to $2D_3+E_3$ and $D_1+S+(n-1)F+E_1$.

Thus, to prove the claim, we may invoke the Nakai-Moishezon Criterion
for ampleness so that it suffices to show that every positive linear
combination of $F$, $D_1$, $D_2$, and $D_3$ has positive self-intersection.
These four divisors have intersection numbers as follows:
\[\begin{array}{r|cccc}
& F & D_1 & D_2 & D_3 \\ \hline
F & 0 & 1 & 2 & 1 \\
D_1 & 1 & n & 2n & n \\
D_2 & 2 & 2n & 4n-2 & 2n-1 \\
D_3 & 1 & n & 2n-1 & n-1 
\end{array}\]
Thus, provided that we choose $n\geq 1$ -- which we may do without loss
of generality -- the claim is proven.  \qed

\vspace{.1in}

We now prove Conjecture~\ref{ratcurve} in the case that the reducible
fibre has a multiple component.  First, assume that $P$ does not lie
on $S$ or the reducible fibre.  The component $C$ of $F$ through $P$
is, by Theorem~\ref{projlines}, a curve of best $D_1$- and
$D_3$-approximation to $P$ (since these divisors contract no curves
through $P$), and it's clear that $C$ is a curve of best
$F$-approximation to $P$.  Finally, if we write $D_2=2D_3+E_3$, then
Corollary~\ref{addeffective} implies that $C$ is a curve of best
$D_2$-approximation to $P$.  Thus, by Corollary~\ref{divsum}, we see that
$C$ is a curve of best $D$-approximation to $P$ for all ample divisors $D$.

Next, assume that $P$ lies on $S$ but not any $E_i$, and consider the
cone $A_S$, generated by:
\[D_1,D_2,D_3,F+D_1,F+D_3,2F+D_2\]
All of these divisors are basepoint free.  The first three contract
$S$, the next two contract no curves through $P$ and map $S$ to a
line, so by Theorem~\ref{projlines}, $S$ is a curve of best
approximation to $P$ with respect to any of these five divisors.
Moreover, $S$ is a curve of best $F$-approximation to $P$, and since
$D_2.S=0$ and $D_2$ is basepoint free, it follows from
Corollary~\ref{addeffective} that $S$ is also a curve of best
$(2F+D_2)$-approximation to $P$.  Thus, by Corollary~\ref{divsum}, it
follows that $S$ is a curve of best $D$-approximation to $P$ for all
$D$ in $A_S$.  If $D$ does not lie in $A_S$, then it must lie in the
cone $B_S$, generated by:
\[F,F+D_1,F+D_3,2F+D_2\]
Since $F$ contracts the component $C$ of $F$ through $P$, we see that
$C$ is a curve of best $F$-approximation to $P$.  Furthermore, $F+D_1$
and $F+D_3$ both contract no curves through $P$ but map $C$ to a line,
so $C$ is a curve of best approximation to $P$ with respect to those
two divisors as well.  Finally, since $C.(2F+D_2)=2=S.(2F+D_2)$, the
fact that $S$ is a curve of best $(2F+D_2)$-approximation to $P$
implies immediately that $C$ is also a curve of best
$(2F+D_2)$-approximation to $P$.  By Corollary~\ref{divsum}, we
conclude that Conjecture~\ref{ratcurve} is true for $P$.

If $P$ lies on exactly one $E_i$ but not $S$, then $F$ and $D_j$ for
$j\neq i$ contract $E_i$, and $D_i$ contracts no curves through $P$
but maps $E_i$ to a line.  Thus, by Theorem~\ref{projlines} and
Corollary~\ref{divsum}, $E_i$ is a curve of best $D$-approximation to $P$
for any ample divisor $D$.

If $P=S\cap E_1$, consider the cone $A_1$, generated by:
\[F+D_1,D_2,D_3,F\]
The last three divisors contract $E_1$, and the first contracts no
curves through $P$ but maps $E_1$ to a line.  We conclude by
Theorem~\ref{projlines} and Corollary~\ref{divsum} that $E_1$ is a
curve of best $D$-approximation to $P$ for any $D$ in $A_1$.  If $D$
does not lie in $A_1$, then it must lie in $B_1$, generated by:
\[D_1,D_2,D_3,F+D_1\]
The first three divisors contract $S$, while the last contracts no
curves through $P$ but maps $S$ to a line.  By Theorem~\ref{projlines}
and Corollary~\ref{divsum}, it follows that $S$ is a curve of best
$D$-approximation to $P$ for any $D$ in $B_1$.  Conjecture~\ref{ratcurve} is
therefore proven for $P$.

If $P=E_1\cap E_2$, then consider the cone $A_2$, generated by:
\[F,D_1,D_3,D_1+D_2\]
The first three contract $E_2$, and the last contracts no curves
through $P$ but maps $E_2$ to a line.  Theorem~\ref{projlines} and
Corollary~\ref{divsum} now imply that $E_2$ is a curve of best
$D$-approximation to $P$ for any $D$ in $A_2$.  If $D$ does not lie in
$A_2$, then it must lie in $B_2$, generated by:
\[F,D_2,D_3,D_1+D_2\]
The first three contract $E_1$, while the last contracts no curves
through $P$ but maps $E_1$ to a line.  Theorem~\ref{projlines} and
Corollary~\ref{divsum} now imply that $E_1$ is a curve of best
$D$-approximation to $P$ for any $D$ in $B_2$.
Conjecture~\ref{ratcurve} is therefore proven for $P$.

Finally, assume that $P=E_2\cap E_3$.  Assume that  $D$ lies in the cone $A_3$,
generated by:
\[F,D_1,D_2,D_2+D_3\]
The first three of these contract $E_3$, and the last contracts no
curves through $P$ but maps $E_3$ to a line.  By
Theorem~\ref{projlines} and Corollary~\ref{divsum}, it follows that
$E_3$ is a curve of best $D$-approximation to $P$ for any $D$ in
$A_3$.  If $D$ does not lie in $A_3$, then it must lie in $B_3$, generated
by:
\[F,D_1,D_3,D_2+D_3\]
The first three contract $E_2$, and the last contracts no curves
through $P$ but maps $E_2$ to a line.  By Theorem~\ref{projlines} and
Corollary~\ref{divsum}, we conclude that $E_2$ is a curve of best
$D$-approximation to $P$ for any $D$ in $B_3$.
Conjecture~\ref{ratcurve} is therefore proven for $X$, for Case~(3),
and for any split rational surface of Picard rank at most four.  \qed

\vspace{.1in}

When the Picard rank is larger than four, the number of different
cases to consider becomes much larger, so we will prove only a few
cases, hopefully representative of the general flavour.  We first turn
our attention to the blowup of $\P^2$ at four $k$-rational points in
general position.  The surface $X$ admits a morphism $\pi_1\colon
X\rightarrow \P^1$ whose fibres are the strict transforms of the
conics through the four blown up points, so that in particular
$\pi_1^*\O(1) = 2L-E_1-E_2-E_3-E_4$.  It also admits a morphism
$\pi_2\colon X\rightarrow \P^2$, which is the blowing down map.  Let
$E_1$, $E_2$, $E_3$, and $E_4$ be the four exceptional divisors of
$\pi_2$, and let $L=\pi_2^*\O(1)$ be the pullback of a line in $\P^2$.

\begin{thm}
Let $X$ be the blowup of $\P^2$ at four $k$-rational points in general
position, and let $P\in X(k)$ be any $k$-rational point.  Then a
sequence of best approximation to $P$ can be chosen to lie as a subset
of the rational curve of minimal degree through $P$.
\end{thm}

\noindent
{\it Proof:} \/ The geometry of $X$ is well understood (see for
instance Example 2.1.2 of [Tsch], or \cite{Ma}).  The nef cone of $X$
is generated by the ten classes $L$, $L_i=L-E_i$,
$D=2L-E_1-E_2-E_3-E_4$, and $D_i=D+E_i$ for $i=1,2,3,4$.  Each of
these is basepoint free: $D$ and each $L_i$ is the fibre of a morphism
to $\P^1$, while $L$ and each $D_i$ corresponds to a birational
morphism to $\P^2$, each one blowing down four pairwise disjoint
smooth rational curves.  There are exactly ten $(-1)$-curves on $X$,
namely $E_i$ and $L-E_i-E_j$ for $i=1,2,3,4$ and $i\neq j$.

We have the following table of
intersection numbers:
\[\begin{array}{r|llllllllll}
  & L & L_1 & L_2 & L_3 & L_4 & D_1 & D_2 & D_3 & D_4 & D \\ \hline
L & 1 & 1 & 1 & 1 & 1 & 2 & 2 & 2 & 2 & 2 \\
L_1 & 1 & 0 & 1 & 1 & 1 & 2 & 1 & 1 & 1 & 1 \\
L_2 & 1 & 1 & 0 & 1 & 1 & 1 & 2 & 1 & 1 & 1 \\
L_3 & 1 & 1 & 1 & 0 & 1 & 1 & 1 & 2 & 1 & 1 \\
L_4 & 1 & 1 & 1 & 1 & 0 & 1 & 1 & 1 & 2 & 1 \\
D_1 & 2 & 2 & 1 & 1 & 1 & 1 & 2 & 2 & 2 & 1 \\
D_2 & 2 & 1 & 2 & 1 & 1 & 2 & 1 & 2 & 2 & 1 \\
D_3 & 2 & 1 & 1 & 2 & 1 & 2 & 2 & 1 & 2 & 1 \\
D_4 & 2 & 1 & 1 & 1 & 2 & 2 & 2 & 2 & 1 & 1 \\
D & 2 & 1 & 1 & 1 & 1 & 1 & 1 & 1 & 1 & 0
\end{array}\]

Let $P\in X(k)$ be any $k$-rational point which does not lie on a
$(-1)$-curve.  Define $B$ to be the cone of nef divisors $A$ such that
$A.D\leq A.L_i$ for all $i$.  Similarly, for each $i$ between 1 and 4,
define the cone $B_i$ as the cone of all nef divisors $A$ such that
$A.L_i\leq A.D$ and $A.L_i\leq A.L_j$ for $j\neq i$.  A short
calculation shows that $B$ is generated by the ten classes:
\[D, D+L, D_i, D+L_i\]
for $i=1,2,3,4$, and that $B_i$ is generated by the ten classes:
\[L, D+L, L_i, L_i+L_j, D+L_i, D_j\]
for $j\neq i$.

Consider the cone $B_i$.  Let $C_i$ be the component of $L_i$ through
$P$, and let $C_D$ be the component of $D$ through $P$.  Each
generator of $B_i$ except $L+D$ either contracts $C_i$, or else it
contracts no curves through $P$ and maps $C_i$ to a line.  Thus, by
Theorem~\ref{projlines}, $C_i$ is a curve of best approximation to $P$
with respect to all the generators of $B_i$ except possibly $L+D$.
For $L+D$, note that $C_i.L=C_i.D=1$, so that a sequence of best $L$-
or $D$-approximation to $P$ along $C_i$ has constant of approximation
equal to 1.  Since $L$ and $D$ are basepoint free, it follows from
Theorem~\ref{product} and Theorem~\ref{projlines} that any sequence
with a better constant of $(L+D)$-approximation to $P$ than $1+1=2$
must have all but finitely many points contained in a union of curves
through $P$ which are contracted by $L$ or $D$.  The only such curve
is $C_D$, which is contracted by $D$ but has $L$-degree 2, and so a
sequence of best $(L+D)$-approximation to $P$ along $C_D$ also has
constant of approximation equal to 2.  Thus, we conclude that both
$C_i$ and $C_D$ are curves of best $(L+D)$-approximation to $P$.  By
Corollary~\ref{divsum}, we conclude that Conjecture~\ref{ratcurve}
is true for any divisor $A$ in $B_i$.

Now consider the cone $B$.  Each generator of $B$ except $L+D$ either
contracts $C_D$, or else contracts no curves through $P$ and maps
$C_D$ to a line.  Thus, by Theorem~\ref{projlines}, $C_D$ is a curve
of best approximation to $P$ with respect to all the generators of
$B$, except possibly $L+D$.  We have already shown that $C_D$ is a
curve of best $(L+D)$-approximation to $P$, so by
Corollary~\ref{divsum}, it follows that Conjecture~\ref{ratcurve} is
true for any ample divisor $A$ in $B$, and thus for any ample divisor
$A$.

There remains the possibility that $P$ might lie on some $(-1)$-curve.
Assume first that $P$ lies on exactly one $(-1)$-curve.  By applying a
suitable automorphism of $X$, we may assume without loss of generality
that this curve is $E_i$.  It is easy to check that any generator of
the cone $B_j$ ($j\neq i$) either contracts $E_i$ or else contracts no
curves through $P$ and maps $E_i$ to a line.  Thus, by
Theorem~\ref{projlines} and Corollary~\ref{divsum}, it follows that
$E_i$ is a curve of best $A$-approximation to $P$ for any $A$ in $B_j$
for $j\neq i$.

If $A$ lies in the cone $B_i$, then consider the subcone $M_i$ of $B_i$
generated by:
\[L+L_i,L_i,L_i+L_j,D_j,D+L_i\]
Each generator is basepoint free, and either contracts $C_i$ or else 
contracts no curves through $P$ and maps $C_i$ to a line.  Thus, by 
Theorem~\ref{projlines} and Corollary~\ref{divsum}, we conclude that
$C_i$ is a curve of best $A$-approximation to $P$ for any $A$ in $M_i$.
If $A$ lies in $B_i$ but not $M_i$, then it must lie in $N_i$, generated by:
\[L,L+L_i,L_i+L_j,D_j,D+L\]
Each generator is basepoint free, and either contracts $E_i$ or else 
contracts no curves through $P$ and maps $E_i$ to a line.  Thus, by 
Theorem~\ref{projlines} and Corollary~\ref{divsum}, we conclude that
$E_i$ is a curve of best $A$-approximation to $P$ for any $A$ in $N_i$,
and therefore for any $A$ in $B_i$.

The only other possibility is that $A$ lies in the cone $B$.  Consider 
the cone $J_i$, generated by:
\[D_i,D_j,L+D,L_j+D,D_i+D\]
where $j\neq i$.  All divisors are basepoint free, and either contract
$E_i$ or else contract no curves through $P$ and map $E_i$ to a point.
By Theorem~\ref{projlines} and Corollary~\ref{divsum}, we conclude
that $E_i$ is a curve of best $A$-approximation to $P$ if $A$ lies in
$J_i$.  If $A$ lies in $B$ but not in $J_i$, then $A$ is in the cone
$J'_i$, generated by:
\[D,D_j,L_i+D,L_j+D,D_i+D\]
where $j\neq i$.  All divisors are basepoint free, and either contract
$D$ or else contract no curves through $P$ and map $D$ to a point.
By Theorem~\ref{projlines} and Corollary~\ref{divsum}, we conclude
that $D$ is a curve of best $A$-approximation to $P$ if $A$ lies in
$J'_i$.  This concludes the proof of Conjecture~\ref{ratcurve} for $P$
lying on a single $(-1)$-curve.

Finally, assume that $P=E_i\cap L-E_i-E_j$ for some $i$ and $j$.  By
symmetry, we may assume that $i=1$ and $j=2$.  Let $A$ be any element
of the nef cone.  Consider the cone $M_1$, generated by:
\[L,L+L_1,L_2,L_3,L_4,L_1+L_3,L_1+L_4,D_1,D_2,D_3+L_4,D_4+L_3,\]
\[D+L_3,D+L_4,D+D_1\]
Each generator of $M_1$ is basepoint free, and either contracts
$E_1$ or else contracts no curves through $P$ and maps
$E_1$ to a line.  Theorem~\ref{projlines} and Corollary~\ref{divsum}
now imply that $E_1$ is a curve of best $A$-approximation to $P$
for any $A$ in $M_1$.  If $A$ does not lie in $M_1$, then it must
lie in the cone $N_1$, generated by:
\[L_1,L_2,L+L_1,L_1+L_3,L_1+L_4,D_2,D_3,D_4,D,D_3+L_4,D_4+L_3,D+L_3,\]
\[D+L_4,D+D_1\]
Each generator of $N_1$ is basepoint free, and either contracts
$L-E_1-E_2$ or else contracts no curves through $P$ and maps
$L-E_1-E_2$ to a line.  By Theorem~\ref{projlines} and
Corollary~\ref{divsum}, we conclude that $L-E_1-E_2$ is a curve of best
$A$-approximation to $P$ for any $A$ in $N_1$, and hence for any ample
divisor $A$.  

This concludes the proof of Conjecture~\ref{ratcurve} for $X$.  \qed

\vspace{.1in}

\begin{thm}
Let $X$ be the blowup of $\P^2$ at five $k$-rational points in general
position.  Conjecture~\ref{ratcurve} is true for $X$.
\end{thm}

\noindent
{\it Proof:} \/ The geometry of $X$ is well understood (see, for
example, \cite{Ma} or \cite{Dr}).  The surface $X$ admits a
blowing-down map $\pi\colon X\rightarrow\P^2$.  Let $L$ be the class
of the preimage of a line, and let $E_1$, $E_2$, $E_3$, $E_4$, and
$E_5$ be the exceptional curves of $\pi$.  The classes $L$ and $E_i$
generate the N\'eron-Severi group of $X$.  The effective cone of $X$
is generated by the 16 $(-1)$-curves on $X$, namely, $E_i$,
$L-E_i-E_j$ for $i\neq j$, and $E=2L-E_1-E_2-E_3-E_4-E_5$.  The ample
cone is dual to the effective cone, and is generated by the following
26 divisors:
\[\begin{array}{ccccc}
L, & L_i=L-E_i, & L_{ij}=E+E_i+E_j, & C_i=E+E_i, & B_i=E+L_i 
\end{array}\]
Let us adopt the convention that the notation $L_{ij}$ implies that
$i<j$.  We have the following table of intersection numbers:

\vspace{.1in}

\noindent
\begin{tabular}{l|ccccc}
            & $L$ & $L_i$ & $L_{ij}$ & $C_i$ & $B_i$ \\ 
\hline
$L$         & 1   & 1     & 2  & 2     & 3     \\
$L_k$       & 1   & $1-\delta_{ik}$ & $1+\delta_{ik}+\delta_{jk}$ 
& $1+\delta_{ik}$ & $2-\delta_{ik}$ \\
$L_{k\ell}$ & 2   & $1+\delta_{ik}+\delta_{i\ell}$ 
& $3-\delta_{ik}-\delta_{j\ell}$
& $2-\delta_{ik}-\delta_{i\ell}$ & $2+\delta_{ik}+\delta_{i\ell}$ \\
$C_k$ & 2 & $1+\delta_{ik}$ & $2-\delta_{ik}-\delta_{jk}$
& $1-\delta_{ik}$ & $1+\delta_{ik}$ \\
$B_k$ & 3 & $2-\delta_{ik}$ & $2+\delta_{ik}+\delta_{jk}$
& $1+\delta_{ik}$ & $2-\delta_{ik}$
\end{tabular}

\vspace{.1in}

Let $P$ be any point on $X$.  We first assume that $P$ does not lie
on any of the 16 $(-1)$-curves on $X$.  In that case, we will show
that for any ample divisor $D$, a sequence of best $D$-approximation
to $P$ lies along a curve in the class $|C_i|$ for some $i$, or in 
the class $|L_i|$ for some $i$.

To do this, we introduce 10 subcones of the nef cone of $X$.  For
$i$ between 1 and 5, we define the cone $M_i$ to be the cone of ample
divisors for which the divisor $L_i$ has minimal degree amongst the
classes $L_j$ and $C_j$, and we define the cone $N_i$ to be the cone
of ample divisors for which the divisor $C_i$ has minimal degree 
amongst the classes $L_j$ and $C_j$.  We will show that for every 
divisor in $M_i$, a sequence of best approximation to $P$ lies along
a curve in $|L_i|$, and similarly that for any divisor in $N_i$, a
sequence of best approximation to $P$ lies along a curve in $|C_i|$.

We first consider $M_i$.  For simplicity, we will consider only $M_1$;
the other four cases differ from this one only by a permutation of the
indices.  The cone $M_1$ has 26 generators:
\[L,\,\,L_1,\,\,L_1+L_i,\,\,L_1+C_i,\,\,L_{ij},\,\,B_1,\,\,L+C_i,
\,\,L_1+C_1,\,\,B_1+L_i\] 
where in all cases we assume that $i\neq 1$ and $j\neq 1$.  Of these
generators, only $L_1$ is not big.  Let $F$ be the unique curve
through $P$ satisfying $F\in |L_1|$.  Since we have assumed that $P$
does not lie on a $(-1)$-curve, it follows that if a curve through $P$
has intersection zero with one of these generators, then the curve in
question must be $F$.

It suffices to show that for any generator of $M_1$, a sequence of
best approximation to $P$ lies on $F$.  For $L$, $L_1$, $L_1+L_i$,
$L_1+C_i$, $L_{ij}$, and $B_1$, this follows from
Theorem~\ref{projlines} because $F$ has degree one or zero.  In all
other cases, $F$ has degree two, and we can write them as follows:
\begin{eqnarray*}
L+C_2 & = & L_3+L_{23} \\
L+C_3 & = & L_2+L_{23} \\
L+C_4 & = & L_2+L_{24} \\
L+C_5 & = & L_2+L_{25} \\
L_1+C_1 & = & L_2+C_2 \\
B_1+L_2 & = & C_3+L_{45} \\
B_1+L_3 & = & C_2+L_{45} \\
B_1+L_4 & = & C_2+L_{35} \\
B_1+L_5 & = & C_2+L_{34}
\end{eqnarray*}
For each generator $D$ in the above list, we have written $D=D_1+D_2$,
where $D_i.F=1$.  Thus, by Theorem~\ref{product}, $F$ contains a sequence
of best $D$-approximation to $P$ unless $P$ lies on some curve of smaller
$D$-degree.  Since we have assumed that $P$ does not lie on a $(-1)$-curve,
and since $D$ lies in $M_1$, it follows that $P$ cannot lie on any curve
of smaller $D$-degree than $L_1$.  (Note that the only candidates for a
curve of smaller degree than two are those which intersect one of $D_1$ or
$D_2$ trivially, and such curves are either $(-1)$-curves or else of the
form $L_i$ or $C_i$ for some $i$.)  We conclude that a sequence of best
$D$-approximation to $P$ lies on $F$.

Now to $N_i$.  Again, we will consider only $N_1$, the cone of
divisors with respect to which $C_1$ has minimal degree amongst the
$L_j$ and $C_j$.  Like $M_1$, $N_1$ has 26 generators:
\[C_1,\,\,C_1+L_i,\,\,C_1+L_{1i},\,\,L_{1i},\,\,B_i,\,\,L+C_1,\,\,
L_1+C_1,\,\,C_1+L_{ij},\,\,B_2+C_2\] 
where again we assume that $i\neq 1$ and $j\neq 1$.  Let $F$ be the
curve in $|C_1|$ containing $P$.  As before, it suffices to show that
for any generator of $N_1$, a sequence of best approximation to $P$
lies on $F$.  For $C_1$, $C_1+L_i$, $C_1+L_{1i}$, $L_{1i}$, and $B_i$,
this follows from Theorem~\ref{projlines} because $F$ has degree one
or zero, and all the generators contract no curves through $P$ other
than $F$.  For the remaining nine generators,
we can write:
\begin{eqnarray*}
L+C_1 & = & L_2+L_{12} \\
L_1+C_1 & = & L_2+C_2 \\
B_2+C_2 & = & B_2+C_2 \\
C_1+L_{45} & = & L_{15}+C_4 \\
C_1+L_{35} & = & L_{15}+C_3 \\
C_1+L_{34} & = & L_{14}+C_3 \\
C_1+L_{25} & = & L_{15}+C_2 \\
C_1+L_{24} & = & L_{14}+C_2 \\
C_1+L_{23} & = & L_{13}+C_2
\end{eqnarray*}
In each case, we have written the generator $D$ as $D=D_1+D_2$, where
$D_i.F=1$.  Thus, by Theorem~\ref{product}, we conclude that since $P$
does not lie on a $(-1)$-curve, a sequence of best $D$-approximation
to $P$ lies on $F$ for any generator $D$ of $N_1$, and by
Corollary~\ref{divsum} also for any $D$ in $N_1$.

It remains only to treat the case where $P$ lies on a $(-1)$-curve.
Since no three $(-1)$-curves meet at a point, $P$ lies either on
exactly one or exactly two $(-1)$-curves.  Let's first assume that $P$
lies on exactly one $(-1)$-curve.  Since the automorphism group of $X$
acts transitively on the set of $(-1)$-curves (see for example
Theorem~2.1 of \cite{Ho}), we may assume that this curve is $E$.

Consider the cone $M_i$.  We will subdivide $M_i$ into two subcones.
The first is $R_i$, generated by the following twenty-two divisor classes:
\[L,L_i,L_i+L_j,L_{j\ell},L_i+C_j,L_i+B_i,L+C_j,L+L_i+C_i\]
where $j$ and $\ell$ are assumed to be different from $i$.  The
classes $L$, $L_i$, $L_i+L_j$, $L_{j\ell}$, $L+C_j$, and $L_i+B_i$ all
either contract the component $F_i$ of $L_i$ through $P$, or else
contract no curves through $P$ and map $F_i$ to a line.  Thus, by
Theorem~\ref{projlines} $F_i$ is a curve of best approximation to $P$
with respect to any of these divisor classes.  For $L+C_j$, we see
that $L.L_i=C_j.L_i=1$, and so Theorem~\ref{product} implies that if
$F_i$ is not a curve of best $(L+C_j)$-approximation to $P$, then a
sequence of best approximation to $P$ can be found on a curve
contracted by either $L$ or $C_j$.  The only such curve through $P$ is
$E$, and we check that $E.(L+C_j)=2$, so $F_i$ is a curve of best
$(L+C_j)$-approximation to $P$ (as is $E$).  This only leaves
$L+L_i+C_i$, which can be rewritten $L+L_i+C_i=L+L_j+C_j$ for any $j$.
Since $L.L_i=L_j.L_i=C_j.L_i=1$ for any $j\neq i$, we may apply
Theorem~\ref{product} to argue that if $F_i$ is not a curve of best
$(L+L_i+C_i)$-approximation to $P$, then by Theorem~\ref{product} a
sequence of best approximation to $P$ must lie on some curve
contracted by $L$, $L_j$, or $C_j$.  The only such curves through $P$
are $E$ and the component of $L_j$ through $P$.  None of these have
$(L+L_i+C_i)$-degree less than $3$ (which is the degree of $F_i$), so
we conclude that $F_i$ is a curve of best $(L+L_i+C_i)$-approximation
to $P$.  Thus, by Corollary~\ref{divsum}, it follows that $F_i$ is a
curve of best $A$-approximation to $P$ for any $A$ in $R_i$.

If $A$ lies in $M_i$ but not $R_i$, then it must lie in the cone $R'_i$,
generated by:
\[L_{j\ell},L_i+C_i,L_i+C_j,B_i,B_i+L_j,B_i+L_i,L+L_i+C_i,L+C_j\]
where $j$ and $\ell$ range over every index different from $i$.  All
but the last two divisors in this list either contract $E$, or else
contract no curves through $E$ and map $P$ to a line.  For the last
two, we see that $E$ has the same degree as $F_i$, and we already know
that $F_i$ is a curve of best $(L+L_i+C_i)$- and
$(L+C_j)$-approximation to $P$.  Thus, we conclude that $E$ is a curve
of best $A$-approximation to $P$ for every $A$ in $R'_i$.

Consider the cone $N_i$.  It is straightforward to check that except
for $L+C_i$, every generator of $N_i$ either contracts $E$ or else
contracts no curves through $P$ and maps $E$ to a line.  We can
rewrite $L+C_i=L_j+L_{ij}$ for any $j\neq i$, and note that
$L_j.E=L_{ij}.E=1$, so that by Theorem~\ref{product}, $E$ is a curve
of best $(L+C_i)$-approximation to $P$ unless $P$ lies on some curve
of smaller $(L+C_i)$-degree which is contracted by $L$ or $C_i$.
Since this is not the case, we conclude by Theorem~\ref{projlines} and
Corollary~\ref{divsum} that $E$ is a curve of best $A$-approximation
to $P$ for any $A$ in $N_i$.

To finish the proof of Conjecture~\ref{ratcurve} for $X$, we just need
to consider the case that $P$ is the intersection point of two
$(-1)$-curves.  Up to automorphisms of $X$, we may assume that one of
these $(-1)$-curves is $E$.  Since the only $(-1)$-curves intersecting 
$E$ are the $E_i$, we may assume that $P=E\cap E_i$ for some $i$.

Consider the cone $J_1$, generated by the following divisor classes:
\[L+L_i,L_i,L_i+L_j,L_{j\ell},L_i+C_j,L_i+B_i\]
where $j$ and $\ell$ range over all values different from $i$.  Each of
these divisors either contracts the component $F_i$ of $L_i$ through $P$,
or else contracts no curves through $P$ and maps $F_i$ to a line.  By
Theorem~\ref{projlines} and Corollary~\ref{divsum}, we conclude that
Conjecture~\ref{ratcurve} is true for any $A$ in $J_1$.

Next, consider the cone $J_2$, generated by:
\[L,L_j,L+L_i,C_i,L_i+L_j,L_{ij},L_{ji},L_{j\ell},L_j+C_\ell,L_{ij}+C_j,
L_j+B_j\] where $j$ and $\ell$ range over all indices different from
$i$ and each other.  Thus, the cone $J_2$ is generated by 41 divisor
classes.  Each of these divisors either contracts $E_i$, or else
contracts no curves through $P$ and maps $E_i$ to a line.  By
Theorem~\ref{projlines} and Corollary~\ref{divsum}, we conclude that
Conjecture~\ref{ratcurve} is true for all $A$ in $J_2$.

Finally, if $A$ does not lie in $J_1$ or $J_2$, then it must lie in $J_3$,
generated by:
\[L_i+B_i,C_i,L_{j\ell},L_i+C_j,L_j+C_\ell,B_j,L_{ij}+C_j,L_j+B_j\]
where again $j$ and $\ell$ are different from each other and $i$, for
a total of 41 divisor classes in the list.  Each class either
contracts $E$, or else contracts no curves through $P$ and maps $E$ to
a line.  Thus, it follows from Theorem~\ref{projlines} and
Corollary~\ref{divsum} that Conjecture~\ref{ratcurve} is true for any
divisor $A$ in $J_3$.

This concludes the proof of Conjecture~\ref{ratcurve} for $X$.  \qed

\section{Further Remarks}

These techniques will probably be able to prove
Conjecture~\ref{ratcurve} for more split rational surfaces, but they
will probably not suffice to prove the conjecture for a general
rational surface.  For example, let $X$ be a blowup of $\P^2$ at six
$k$-rational points in general position, embedded in $\P^3$ as a
smooth cubic surface.  The family $\mathcal{F}$ of plane cubic curves
passing through the set of blown up points in $\P^2$ is a
three-dimensional (projective) linear subspace of the $\P^{9}$ of
plane cubic curves.  The Zariski closure $Z$ of the set of cuspidal
cubics is a closed subset of dimension 7.  Therefore,
$Z\cap\mathcal{F}$ has dimension at least one, so if the six points
are chosen to lie on some cuspidal cubic whose cusp is not one of the
six blown up points (which will generically be the case), there will
be a one-dimensional family $\mathcal{C}$ of cuspidal cubic curves on
$X$.

Let $P$ be the cusp of a cuspidal cubic $C$, and assume that $P$ is
$k$-rational.  Then one can find a sequence on $C$ which approximates
$P$ with constant of approximation $3/2$.  In particular, no rational
curve of minimal degree through $P$ is a curve of best approximation
to $P$, since that curve will generally be a conic ($P$ cannot lie on
any of the 27 lines on $X$).  This suggests that any proof of
Conjecture~\ref{ratcurve} for $X$ will be beyond the techniques of
this paper.

Note also that Vojta's Main Conjecture (\cite{Vo}, Conjecture~3.4.3)
implies Conjecture~\ref{ratcurve} for many varieties.  Vojta's Main
Conjecture is as follows:

\begin{conj}[Vojta's Main Conjecture]\label{vstrong}
Let $X$ be a smooth algebraic variety defined over a number field $k$,
with canonical divisor $K$.  Let $S$ be a finite set of places of $k$.
Let $L$ be a big divisor on $X$, and let $D$ be a normal
crossings divisor on $X$.  Choose height functions $h_K$ and $h_L$ for
$K$ and $L$, respectively, and define a proximity function $m_S(D,P) =
\sum_{v\in S} h_{D,v}(P)$ for $D$ with respect to $S$, where $h_{D,v}$
is a local height function for $D$ at $v$.  Choose any $\epsilon>0$.
Then there exists a nonempty Zariski open set $U=U(\epsilon)\subset X$
such that for every $k$-rational point $P\in U(k)$, we have the
following inequality:
\begin{equation}
M_S(D,P)H_K(P) \leq H_L(P)^{\epsilon}
\end{equation}
\end{conj}

(Note that our notation is multiplicative, rather than the additive
notation used in \cite{Vo}.  Thus, Vojta's $h$ is our $\log H$, and
Vojta's $m$ is our $\log M$.)

Assume that $X$ is as in the conjecture, and has non-negative Kodaira
dimension.  Let $P$ be any fixed $k$-rational point, lying on a
rational curve $C\subset X$, and let $D$ be any normal crossings
divisor containing $P$.  Since $X$ has non-negative Kodaira dimension,
there is a non-empty open set $U_1\subset X$ such that there is a
constant $c$ satisfying $H_K(Q)>c>1$ for all $Q\in U_1(k)$.
Furthermore, if $S$ is the set of archimedean places of $k$, we may
write:
\[M_S(D,Q)\gg\left(\dist(P,Q)\right)^{-1}\]
where the implied constant is independent of $Q$.  (If $Q$ lies on the
support of $D$, then $M_S(D,Q)$ is infinite, but this will not affect
our proof.)

Let $d$ be the constant of approximation for $P$ on $C$.  For any
$\epsilon>0$, let $U_\epsilon$ be the intersection of $U_1$ with the
non-empty open set provided by Conjecture~\ref{vstrong}.  Then for any
point $Q\in U_\epsilon(k)$, Conjecture~\ref{vstrong} implies:
\[\dist(P,Q)H_L(P)^\epsilon>c\]
If we choose $\epsilon<d$, then it's clear that any sequence of best
approximation to $P$ must eventually lie in the complement of
$U_\epsilon$.  If $X$ is a surface, then this complement is a union of
curves, in which case the truth of Conjecture~\ref{ratcurve} is clear.
If $X$ is of higher dimension, then a simple Noetherian induction
shows that it suffices to assume that $X$ contains no subvariety $V$
whose Kodaira dimension is negative but whose dimension is at least
two.  Slightly more generally, we have therefore proven the following
theorem:

\begin{thm}
Let $V$ be a smooth algebraic variety with non-negative Kodaira
dimension.  Let $P\in V(k)$ be any rational point with a rational
curve $C/k$ on $V$ through $P$, and assume that every subvariety of
$V$ with negative Kodaira dimension satisfies
Conjecture~\ref{ratcurve} for $P$. Assume that
Conjecture~\ref{vstrong} is true for $V$.  Then
Conjecture~\ref{ratcurve} is true for $P$ on $V$.
\end{thm}

One can also use our techniques to prove Conjecture~\ref{ratcurve} for
certain special points on some varieties of nonnegative Kodaira
dimension.  For example, let $X$ be a smooth quartic in $\P^3$
containing a line $L$.  This $X$ is a K3 surface, and $L$ is a smooth
rational curve on $X$ with self-intersection $-2$.  For a general such
$X$, the effective cone of $X$ is spanned by an elliptic curve $E$ and
$L$.  (Indeed, one can choose $E$ so that $E\cup L$ is a hyperplane
section of $X$.  See section~5 of \cite{Ko} for details.)

The closure of the nef cone of $X$ is therefore spanned by the
divisors $E$ and $D=2E+3L$.  The former corresponds to a morphism from
$X$ to $\P^1$ giving $X$ the structure of an elliptic surface, and $D$
corresponds to the contraction of the $-2$ curve $L$.  Let $P$ be a
point on $L$.  Then clearly any sequence of best $D$-approximation to
$P$ lies on $L$, and by Theorem~\ref{ratapprox}, a sequence of best
$(E+L)$-approximation to $P$ also lies along $L$.  Therefore, for any
ample divisor $A$ in the positive span of $D$ and $E+L$, a sequence of
best $A$-approximation to $P$ must lie along $L$, by
Theorem~\ref{divsum}.

\end{document}